
\def\LaTeX{%
  \let\Begin\begin
  \let\End\end
  \let\salta\relax
  \let\finqui\relax
  \let\futuro\relax}

\def\UK{\def\our{our}\let\sz s}
\def\USA{\def\our{or}\let\sz z}

\UK 



\LaTeX

\USA


\salta

\documentclass[twoside,12pt]{article}
\setlength{\textheight}{24cm}
\setlength{\textwidth}{16cm}
\setlength{\oddsidemargin}{2mm}
\setlength{\evensidemargin}{2mm}
\setlength{\topmargin}{-15mm}
\parskip2mm


\usepackage[usenames,dvipsnames]{color}
\usepackage{amsmath}
\usepackage{amsthm}
\usepackage{amssymb,bbm}
\usepackage[mathcal]{euscript}

\usepackage{cite}
\usepackage{hyperref}
\usepackage{enumitem}

\usepackage[ulem=normalem,draft]{changes}
%
%

%
 
\definecolor{ciclamino}{rgb}{0.5,0,0.5}
\definecolor{blu}{rgb}{0,0,0.7}
\definecolor{rosso}{rgb}{0.85,0,0}

\def\juerg #1{{\color{green}#1}}
\def\pier #1{{\color{rosso}#1}} 
\def\an #1{{\color{blue}#1}}
\def\last #1{{\color{rosso}#1}}

\def\juerg #1{#1}
\def\an #1{{#1}}
\def\last #1{{#1}}
\def\pier #1{{#1}}




\bibliographystyle{plain}


%
\newtheorem{theorem}{Theorem}[section]

\newtheorem{definition}[theorem]{Definition}

\newtheorem{lemma}[theorem]{Lemma}

\finqui

\def\Bthm{\Begin{theorem}}
\def\Ethm{\End{theorem}}
\def\Blem{\Begin{lemma}}
\def\Elem{\End{lemma}}
\def\Bprop{\Begin{proposition}}
\def\Eprop{\End{proposition}}
\def\Bcor{\Begin{corollary}}
\def\Ecor{\End{corollary}}
\def\Brem{\Begin{remark}\rm}
\def\Erem{\End{remark}}

\def\Bdim{\Begin{proof}}
\def\Edim{\End{proof}}
\def\Bcenter{\Begin{center}}
\def\Ecenter{\End{center}}
\let\non\nonumber




\def\step #1 \par{\medskip\noindent{\bf #1.}\quad}
\def\jstep #1: \par {\vspace{2mm}\noindent\underline{\sc #1 :}\par\nobreak\vspace{1mm}\noindent}

\def\Lip{Lip\-schitz}

\def\lhs{left-hand side}
\def\rhs{right-hand side}




\def\multibold #1{\def\arg{#1}%
  \ifx\arg\pto \let\next\relax
  \else
  \def\next{\expandafter
    \def\csname #1#1#1\endcsname{{\bf #1}}%
    \multibold}%
  \fi \next}

\def\pto{.}

\def\multical #1{\def\arg{#1}%
  \ifx\arg\pto \let\next\relax
  \else
  \def\next{\expandafter
    \def\csname cal#1\endcsname{{\cal #1}}%
    \multical}%
  \fi \next}


\def\multimathop #1 {\def\arg{#1}%
  \ifx\arg\pto \let\next\relax
  \else
  \def\next{\expandafter
    \def\csname #1\endcsname{\mathop{\rm #1}\nolimits}%
    \multimathop}%
  \fi \next}

\multibold
qwertyuiopasdfghjklzxcvbnmQWERTYUIOPASDFGHJKLZXCVBNM.

\multical
QWERTYUIOPASDFGHJKLZXCVBNM.

\multimathop
diag dist div dom mean meas sign supp .


\def\Accorpa #1#2 #3 {\gdef #1{\eqref{#2}--\eqref{#3}}%
  \wlog{}\wlog{\string #1 -> #2 - #3}\wlog{}}


\def\<#1>{\mathopen\langle #1\mathclose\rangle}
\def\norma #1{\mathopen \| #1\mathclose \|}

\def\I2 #1{\int_{Q_t}|{#1}|^2}
\def\IT2 #1{\int_{Q^t}|{#1}|^2}
\def\IO2 #1{\norma{{#1(t)}}^2}
\def\ov #1{{\overline{#1}}}
\def\next{\\ & \quad}
\def\aeO{\checkmmode{a.e.\ in~$\Omega$}}

\def\iO{\int_\Omega}

\def\dt{\partial_t}
\def\dn{\partial_{\bf n}}
\def\S{{\cal S}}

\def\X{{\cal X}}
\def\Y{{\cal Y}}
\def\Uh{{\cal U}}

\def\checkmmode #1{\relax\ifmmode\hbox{#1}\else{#1}\fi}


\def\bu{{\bf u}}
\def\bh{{\bf h}}

\def\erre{{\mathbb{R}}}

\def\enne{{\mathbb{N}}}

\def\J{{\cal J}}




\def\genspazio #1#2#3#4#5{#1^{#2}(#5,#4;#3)}
\def\spazio #1#2#3{\genspazio {#1}{#2}{#3}T0}

\def\L {\spazio L}
\def\H {\spazio H}

\def\C #1#2{C^{#1}([0,T];#2)}



\def\Lx #1{L^{#1}(\Omega)}
\def\Hx #1{H^{#1}(\Omega)}

\def\Ldue{\Lx 2}

\def\Huno{\Hx 1}
\def\Hdue{\Hx 2}

\def\Liq{{L^\infty(Q)}}




\let\vp\varphi

\def\a{\alpha}	
\def\b{\beta}	
\def\s{\sigma}  
\def\m{\mu}	    
\def\ph{\varphi}	

\def\h{\mathbbm{h}}

\let\TeXchi\chi                         
\newbox\chibox
\setbox0 \hbox{\mathsurround0pt $\TeXchi$}
\setbox\chibox \hbox{\raise\dp0 \box 0 }
\def\chi{\copy\chibox}



\def\ubar{\overline{\bf u}}
\def\uebar{\overline{u}_1}
\def\uzbar{\overline{u}_2}

\def\phial{\varphi_\gamma}
\def\sial{\sigma_\gamma}

\def\phialn{\varphi_{\gamma_n}}
\def\sialn{\sigma_{\gamma_n}}                  
\def\mual{\mu_\gamma}
\def\mualn{\mu_{\gamma_n}}

\def\pal{p_\gamma}
\def\qal{q_\gamma}
\def\paln{p_{\gamma_n}}
\def\qaln{q_{\gamma_n}}
\def\ral{r_\gamma}
\def\raln{r_{\gamma_n}}

\def\CP0{(${\mathcal{CP}}_0$)}
\def\CPal{(${\mathcal{CP}}_\gamma$)}
\def\CPaltil{($\widetilde{{\mathcal{CP}}}_\gamma$)}

\def\SAL{{\cal S}_\gamma}
\def\SO{{\cal S}_0}

\let\hat\widehat

\def\fal{F_{1,\gamma}}
\def\faln{F_{1,\gamma_n}}

\def\uad{{\cal U}_{\rm ad}}
\def\UR{{\cal U}_{R}}
\def\bmu{\overline\mu}
\def\bvp{\overline\varphi}
\def\bsigma{\overline\sigma}

\def\bph{{\ov \ph}}
\def\bm{{\ov \m}}   
\def\bs{{\ov \s}}

\def\indi{I_{[-1,1]}}
\def\sindi{\partial I_{[-1,1]}}

\usepackage{amsmath}
\DeclareFontFamily{U}{mathc}{}
\DeclareFontShape{U}{mathc}{m}{it}%
{<->s*[1.03] mathc10}{}

\DeclareMathAlphabet{\mathscr}{U}{mathc}{m}{it}
\Begin{document}


%
\title{Optimal control problems \an{with sparsity} for phase field 
\linebreak tumor growth models involving variational inequalities}
\author{}
\date{}
\maketitle
\Bcenter
\vskip-1.5cm
{\large\sc Pierluigi Colli$^{(1)}$}\\
{\normalsize e-mail:{\tt pierluigi.colli@unipv.it}}\\[0.25cm]
{\large\sc Andrea Signori$^{(1)}$}\\
{\normalsize e-mail:{\tt andrea.signori02@universitadipavia.it}}\\[0.25cm]
{\large\sc J\"urgen Sprekels$^{(2)}$}\\
{\normalsize e-mail: {\tt juergen.sprekels@wias-berlin.de}}\\[.5cm]
$^{(1)}$
{\small Dipartimento di Matematica ``F. Casorati''}\\
{\small Universit\`a di Pavia}\\
{\small via Ferrata 5, I-27100 Pavia, Italy}\\[.3cm] 
$^{(2)}$
{\small Department of Mathematics}\\
{\small Humboldt-Universit\"at zu Berlin}\\
{\small Unter den Linden 6, D-10099 Berlin, Germany}\\[2mm]
{\small and}\\[2mm]
{\small Weierstrass Institute for Applied Analysis and Stochastics}\\
{\small Mohrenstrasse 39, D-10117 Berlin, Germany}\\[10mm]
\Ecenter
\Begin{abstract}
\noindent 
This paper treats a distributed optimal control problem for a tumor growth model of Cahn--Hilliard type including chemotaxis. The evolution 
of the tumor fraction is governed by a variational inequality  
corresponding to a double obstacle nonlinearity occurring in the associated potential. In addition, the control and state variables are nonlinearly coupled \an{and, furthermore,} the cost functional contains a nondifferentiable
term like the $L^1$--norm in order to include sparsity effects
\an{which is of utmost relevance, especially time sparsity, in the context of cancer therapies as applying a control to the system reflects in exposing the patient to \last{an} intensive medical treatment.}
To cope with the difficulties originating from the variational inequality in the state system, we employ the so-called ``deep quench approximation'' in
which the convex part of the double obstacle 
potential is approximated by logarithmic functions. For such functions, first-order necessary conditions of optimality can be established 
by invoking recent results. We use these results to  derive corresponding optimality 
conditions also for the double obstacle case, by deducing a variational inequality in 
terms of the associated adjoint state variables. The resulting variational
inequality can be exploited to also \pier{obtain} sparsity results for the optimal 
controls. 

\End{abstract}
\salta
\pagestyle{myheadings}
\newcommand\testopari{\sc Colli -- Signori -- Sprekels}
\newcommand\testodispari{\sc Optimal control of a variational inequality tumor growth model}
\markboth{\testopari}{\testodispari}
\finqui
%

\section{Introduction}
\label{INTRO}
\setcounter{equation}{0}

Let  $\a>0,~\b>0$, and let $\Omega\subset\erre^3$ denote some open and bounded domain having a smooth boundary $\Gamma=\partial\Omega$ and the unit outward normal $\,{\bf n}$.
We denote by $\dn$ the outward normal derivative to $\Gamma$.  Moreover, we fix some final time $T>0$ and
introduce for every $t\in (0,T)$ the sets $Q_t:=\Omega\times (0,t)$ and $Q^t:=\Omega\times (t,T)$.
\last{Furthermore}, we set $Q:=Q_T$ and $\Sigma:=\Gamma\times (0,T)$.
We then consider the following optimal control problem: 

\vspace{3mm}\noindent
(${\cal CP}$) \quad Minimize the cost functional
\begin{align} 
	&	\J((\mu,\vp,\sigma),{\bf u})
   :=\,  \frac{b_1}2 \iint_Q |\vp-\widehat \vp_Q|^2
  + \frac{b_2}2 \iO |\vp(T)-\widehat\vp_\Omega|^2
  + \frac{b_0}2 \iint_Q |{\bf u}|^2 
  \,+\,\kappa\,g({\bf u})\nonumber\\
     & =:\, \J_1((\mu,\vp,\sigma),{\bf u}) + \kappa g(\bu)
	\label{cost}
\end{align} 
subject to  the state system 
\begin{align}
\label{ss1}
&\alpha\dt\mu+\dt\ph-\Delta\mu=P(\ph)(\sigma+\chi(1-\ph)-\mu) - \h(\ph)u_1 \quad&&\mbox{in }\,Q\,,\\[1mm]
\label{ss2}
&\beta\dt\vp-\Delta\vp+F_1'(\vp) + F_2'(\vp)=\mu+\chi\,\sigma \quad
&&\mbox{in }\,Q\,,\\[1mm]
\label{ss3}
&\dt\sigma-\Delta\sigma=-\chi\Delta\vp-P(\ph)(\sigma+\chi(1-\ph)-\mu)+u_2\quad&&\mbox{in }\,Q\,,\\[1mm]
\label{ss4}
&\dn \mu=\dn\vp=\dn\sigma=0 \quad&&\mbox{on }\,\Sigma\,,\\[1mm] 
\label{ss5}
&\mu(0)=\mu_0,\quad \vp(0)=\vp_0,\quad \sigma(0)=\sigma_0\,\quad &&\mbox{in }\,\Omega\,,
\end{align}
\Accorpa\Statesys {ss1} {ss5}
and to the control constraint
\begin{equation}
\label{cont:constr}
{\bf u}=(u_1,u_2)\in\uad\,.
\end{equation}
Here, $b_1,b_2, \kappa$ are nonnegative constants, while $\,b_0\,$  
is positive. Moreover,
 $\widehat\vp_Q$ and $\widehat \vp_\Omega$ are given target functions, and $g$ denotes a convex 
but not necessarily differentiable  functional that may account for possible sparsity effects;
a typical case is $g(\bu)=\|\bu\|_\last{{L^1(Q)}^2}$. Moreover, the set of admissible controls $\uad$ is a nonempty, closed and convex subset of the control space
\begin{equation}
\label{defU}
\Uh:= L^\infty(Q)^2.
\end{equation}

The state system \Statesys\ constitutes a simplified and relaxed version of the four-species thermodynamically consistent model for tumor growth originally proposed by Hawkins-Daruud 
et al.\ in \cite{HZO} that additionally includes chemotaxis effects. 
Let us briefly review the role of the occurring symbols. The primary variables $\ph, \m, \s$ 
denote the tumor fraction, the associated chemical potential, and the nutrient concentration, respectively.
Furthermore, the additional term $\a\dt\m$ corresponds to a parabolic regularization of equation \eqref{ss1},
while $\b\dt\ph$ is the viscosity contribution to the Cahn--Hilliard equation.
The nonlinearity $P$ denotes a proliferation function, whereas the positive constant $\chi$
represents the chemotactic sensitivity
\last{and provides the system with a cross-diffusion coupling.}
The evolution of the tumor fraction is mainly governed by
the nonlinearities $F_1$ and $F_2$ \last{whose derivatives occur} in \eqref{ss2}. Here, $F_2$ is smooth, typically a 
concave function. As far as
$F_1$ is concerned, we consider in this paper the functions   
\begin{align}
\label{logpot} 
&F_{\rm 1,log}(r)=\left\{\begin{array}{ll}
(1+r)\,\ln(1+r)+(1-r)\,\ln(1-r)\quad&\mbox{for $\,r\in (-1,1)$}\\
2\,\ln(2)\quad&\mbox{for $\,r\in\{-1,1\}$\,,}\\
+\infty\quad&\mbox{for $\,r\not\in[-1,1]$}
\end{array}\right.\\[1mm]
&\indi (r)=\left\{\begin{array}{ll}
0\quad&\mbox{for $\,r\in [-1,1]$}\\
+\infty\quad&\mbox{for $\,r\not\in[-1,1]$}
\end{array}\right. \,.
\label{obspot}
\end{align}
We assume that $\indi+F_2$ is a double-well potential. This is \pier{actually} the case 
if $F_2(r)=k(1-r^2)$, where $k>0$; the function
$\indi+F_2$ is then referred to as a 
{\em double obstacle} potential. Note also that $F'_{\rm 1,log}(r)$ becomes 
unbounded as $r\searrow -1$ and $\an{r\nearrow 1}$,
and that in the case of \eqref{obspot} the second equation \eqref{ss2} has to be interpreted
 as a differential inclusion,
where $F_1'(\ph)$ {is} understood in  the sense of \an{subdifferentials. Namely, \eqref{ss2} has to be written
as}
\begin{equation}
\label{ss2new}
\beta\dt\ph-\Delta\ph+\xi+F_2'(\vp)=\mu+\chi\sigma, \quad \xi\in\sindi(\ph).
\end{equation}

The control variable $u_1$, which is nonlinearly coupled to the state variable $\varphi$
in the phase equation \eqref{ss1}, models the application of a cytotoxic drug to 
the system; it is multiplied by a truncation function $\h({\cdot})$ in order to have the action
only in the spatial region where the tumor cells are located. Typically,
one assumes that $\h(-1)=0, \linebreak \h(1)=1$, \an{and $\h(\ph)$} is in between if $-1<\ph<1$;
see \cite{GLSS, GARL_1, HKNZ, KL} for some insights on possible choices of $\h$.
On the other hand, the control $u_2$ can model  
either an external medication or some nutrient supply.

As far as well-posedness {is concerned}, the above model was already investigated in the 
case $\chi=0$ in \cite{CGH,CGRS1,CGRS2,CRW},
and in \cite{FGR} with $\a=\b=\chi=0$.
There the authors also pointed out how \last{the relaxation parameters} $\a$ and $\b$ can be set to zero, by providing
the proper framework in which a limit system can be identified and uniquely solved.
We also note that in \cite{CGS24} a version has been studied in which the Laplacian in the equations 
\eqref{ss1}--\eqref{ss3} has been replaced by fractional powers of a more general class of selfadjoint operators
having compact resolvents. A model which is similar to the one studied in this note was the subject
of \cite{CSS1,ST}.

For some nonlocal variations of the above model we refer to \cite{FLRS, FLS, SS}.
Moreover, in order to better emulate in-vivo tumor growth,
it is possible to include in similar models the effects generated by the fluid flow development
by postulating a Darcy{'s} law or a Stokes--Brinkman{'s} law.
In this direction, we refer to \cite{WLFC,GLSS,DFRGM ,GARL_1,GARL_4,GAR, EGAR, FLRS,GARL_2, GARL_3}, and {we also mention} \cite{GLS}, where elastic effects are included.
For further models, discussing the case of multispecies,
we refer the reader to \cite{DFRGM,FLRS}.

The investigation of associated optimal control problems also presents a 
wide number of results of which we mention{\cite{CGRS3,EK, EK_ADV, S, S_a, S_b, S_DQ, SigTime, FLS, ST, CGS24, CSS1, GARLR, KL, SW, GLS_OPT}. 
The optimal control problem (${\cal CP}$) has recently been investigated by the present authors in \cite{CSS2} for the
case of regular or logarithmic nonlinearities $F_1$. For such nonlinearities, well-posedness of the state system
\Statesys, suitable differentiability properties of the control-to-state mapping, the existence of optimal controls, 
as well as first-order necessary and second-order sufficient optimality conditions could be established. 
In this paper, we focus on the nondifferentiable case when $F_1=\indi$. While a well-posedness result 
was proved in \cite{CSS2} also for
this case (in which \eqref{ss2} has to be replaced by the inclusion \eqref{ss2new}), the corresponding optimal 
control problem has not yet been treated. While the existence of optimal controls is not too difficult to show, the  
derivation of necessary optimality is challenging since standard constraint qualifications to establish the 
existence of suitable Lagrange multipliers are not available. In order to handle this difficulty, we 
employ the so-called ``deep quench approximation'' which has proven to be a useful tool in a number 
of optimal control problems for Cahn--Hilliard systems involving double obstacle potentials (cf., e.g., the papers 
\cite{CFGS, CGSEECT,CGSconv, CGS24,CGSCC, CGS2021,S_DQ}). 

In all of these works, the starting point was that the optimal control problem (we will later
denote this problem by \CPal) had been successfully treated (by proving Fr\'echet differentiability of 
the control-to-state operator and establishing
first-order necessary optimality conditions in terms of a variational inequality and the adjoint
state system) 
for the case when in the state system \Statesys\ the nonlinearity 
$F_1$ \last{is, for $\gamma>0$,} given by 
\begin{equation}
F_{1,\gamma}:=\gamma\, F_{\rm 1,log}.
\label{F1G}
\end{equation}
We obviously have that
\begin{align}
\label{monoal}
&0\,\le\,F_{1,{\gamma_1}}(r)\,\le\,F_{1,{\gamma_2}}(r)\quad\forall\,r\in\erre, \quad\mbox{if }\,0<\gamma_1<\gamma_2, 
\\[1mm]
\label{limF1al}
&\lim_{\gamma\searrow0} F_{1,\gamma}(r)\,=\,\indi(r)\quad\forall\,r\in \erre.
\end{align}
In addition, \pier{we note that} $\,F_{\rm 1,log}'(r)=\ln\left(\frac{1+r}{1-r}\right)$ \,and\, $F_{\rm 1,log}''(r)=\frac 2{1-r^2}>0$\, for
$r\in (-1,1)$, and thus, in particular,
\begin{align}
\label{limF1al1}
&\lim_{\gamma\searrow0}\, F_{1,\gamma}'(r)= \an{\lim_{\gamma\searrow0}\, \gamma\, F_{1,\rm log}'(r) = 0} \quad\mbox{for }\,-1<r<1,\\
\label{limF1al2}
&\lim_{\gamma\searrow 0} \Bigl(\,\lim_{r\searrow -1}\fal'(r)\Bigr)=-\infty, \quad
\lim_{\gamma\searrow0} \Bigl(\,\lim_{\an{r\nearrow 1}}\fal'(r)\Bigr)=+\infty.
\end{align} 
We may therefore regard the graphs of the single-valued functions 
\begin{equation}
F_{1,\gamma}'(r)\,=\, \gamma\, F_{\rm 1,log}'(r), \quad \mbox{for}\quad r\in (-1,1)\quad\mbox{and}\quad \gamma>0,
\end{equation} 
as approximations to the graph of the multi-valued subdifferential $\sindi$ 
from the interior of $(-1,1)$.

For both $F_1=\indi$ (in which case \eqref{ss2} has to be replaced by 
the inclusion \eqref{ss2new}) and $F_1=\fal$ (where $\gamma>0$), the well-posedness 
results from \cite{CSS2} yield the existence of a unique solution
$(\mu,\vp,\sigma)$ and $(\mual,\phial,\sial)$ to the state system
\Statesys\ provided that the controls $\bu=(u_1,u_2)$ 
belong to $L^\infty(0,T;\pier{L^2(\Omega)})^2$. It is  natural to expect that $(\mual,
\phial,\sial)\to (\mu,\vp,\sigma)$ as $\gamma\searrow0$ in a suitable topology.

Below (cf.~Theorem \ref{THM:DQ}), we will show that this is actually true.
Owing to the construction, the approximating functions} $\,\phial\,$  automatically attain their values in
the domain of $\indi$; that is, we have $\|\phial\|_{L^\infty(Q)}\,\le\,1$ \,for all $\gamma>0$. 

Let us now consider the control problem, which in the following will be denoted by \CP0\ if
$F_1=\indi$ and by \CPal\ if $F_1=\fal$. The general strategy is then to derive 
uniform (with respect to $\gamma\in (0,1]$) a priori estimates for the state and adjoint state variables of 
an ``adapted'' version of \CPal\ that are sufficiently strong as to permit a passage to the 
limit as $\gamma\searrow0$ in order to derive meaningful first-order necessary optimality conditions also for \CP0. 
It turns out that this strategy succeeds. 

Another \an{remarkable} novelty of this paper is the discussion of the sparsity of optimal controls for \CP0. 
Since the seminal paper \cite{stadler2009}, sparse optimal controls have been discussed extensively in the literature. Directional sparsity was introduced in \cite{herzog_stadler_wachsmuth2012, herzog_obermeier_wachsmuth2015} and extended to semilinear parabolic optimal control problems in
\cite{casas_herzog_wachsmuth2017}. Sparse optimal controls for reaction-diffusion equations  were investigated in \cite{casas_ryll_troeltzsch2013,casas_ryll_troeltzsch2015}. In the recent work \cite{ST}, sparsity results that apply to nonlinearities $\,F_1\,$ of logarithmic type were
established for a slightly different state system.
In view of the medical background, the focus in \cite{ST} was set on sparsity with 
respect to time, since temporal sparsity means that
the controls (e.g., cytotoxic drugs) are not needed in certain time periods. It turns out that the technique used in 
\cite{ST} can be adapted to establish 
sparsity results also for our state system for the 
nondifferentiable case $F_1=\indi$ in which
the evolution of the tumor fraction is governed by a variational inequality. 
The results obtained, however, are weaker than those recovered in \cite{ST}
for the differentiable case. This is not entirely unexpected in view of the fact
that less information on the adjoint state variables can be recovered from
the corresponding adjoint state system than in the simpler differentiable
situation. 

The paper is organized as follows: in Section 2, we collect auxiliary results  on the state system \Statesys\
that have been established in \cite{CSS2}. The subsequent Section 3 brings a detailed analysis of the
deep quench approximation. Section 4 is then devoted to the derivation of first-order necessary optimality 
conditions for the case $F_1=\indi$. In the final Section 5, we investigate sparsity properties
of the optimal controls for the double obstacle case.   
 
Throughout the paper, we make repeated use of H\"older's inequality, of the elementary \pier{Young} inequality
\begin{equation}
\label{Young}
a b\,\le \delta |a|^2+\frac 1{4\delta}|b|^2\quad\forall\,a,b\in\erre, \quad\forall\,\delta>0,
\end{equation}
as well as the continuity of the embeddings $H^1(\Omega)\subset L^p(\Omega)$ for $1\le p\le 6$ and 
$\Hdue\subset C^0(\overline\Omega)$.


\section{General setting  and properties of the control-to-state operator}
\label{STATE}
\setcounter{equation}{0}

In this section, we introduce the general setting of our control 
problem and state some results on the state system \eqref{ss1}--\eqref{ss5} that in the present form have been
established in \cite{CSS2}. For similar results, we also refer to the papers \cite{CSS1}
and \cite{ST}.
 
To begin with, for a Banach space $\,X\,$ we denote by $\|\cdot\|_X$
the norm in the space $X$ or in a power {thereof}, \an{by $\,X^*\,$ its dual space, and with $\< \, \cdot\, , \, \cdot \, >_{X}$ the duality pairing between \pier{$X^*$ and $X$}. 
For any $1 \leq p \leq \infty$ and $k \geq 0$, we denote the standard Lebesgue and Sobolev spaces on $\Omega$ by $L^p(\Omega)$ and $W^{k,p}(\Omega)$, and the corresponding norms by $\norma{\,\cdot\,}_{L^p(\Omega)}=\norma{\,\cdot\,}_{p}$ and $\norma{\,\cdot\,}_{W^{k,p}(\Omega)}$, respectively. 
For the case $p = 2$, these become Hilbert spaces and we employ the standard notation $H^k(\Omega) := W^{k,2}(\Omega)$. }
As usual, for Banach spaces $\,X\,$ and $\,Y\,$ we introduce the linear space
$\,X\cap Y\,$ which becomes a Banach space when equipped with its natural norm \an{$\,\|v\|_{X \cap Y}:=
\|v\|_X\,+\,\|v\|_Y$, for $\,v\in X\cap Y$.}
Moreover, we recall the definition~\eqref{defU} of \an{the control space} ${\cal U}$ and introduce the spaces
\begin{align}
  & H := \Ldue \,, \quad  
  V := \Huno\,,   \quad
  W_{0} := \{v\in\Hdue: \ \dn v=0 \,\mbox{ on $\,\Gamma$}\}.
  \label{defHVW}
\end{align}
Furthermore, by $\,(\,\cdot\,,\,\cdot\,)$ \pier{and  $\,\Vert\,\cdot\,\Vert$ we denote the standard inner product 
and related norm in $\,H$, and for simplicity  we also set $\<\,\cdot\,,\,\cdot\,> := \<\,\cdot\,,\,\cdot\,>_V $.}

We make the following assumptions on the data of the system.
\begin{enumerate}[label={\bf (A\arabic{*})}, ref={\bf (A\arabic{*})}]
\item \label{const:weak}
	\last{$\alpha,\beta $, and $\chi$} are positive constants.
\item \label{F:weak}
	$F=F_1+F_2$, where $\,F_1:\erre\to [0,+\infty]\,$ is convex and lower semicontinuous with
	$\,F_1(0)=0$, and where $F_2 \in C^5(\erre)$ has a \Lip\ continuous derivative $F'_2$.
\item \label{P:weak}
	$P ,\h\in C^3(\erre)\cap W^{3,\infty}(\erre)$ are nonnegative and bounded.
\item\label{uad}
	With fixed given constants $\underline u_i,\hat u_i$ satisfying $\underline u_i<\hat u_i$, $i=1,2$, we have	
\begin{equation}
\label{defuad}
\uad=\left\{\bu=(u_1,u_2)\in \Liq^2: \underline u_i\le u_i\le\hat u_i\,\mbox{ a.e. in }\,Q\,\mbox{ for }\,i=1,2\right\}. 
\end{equation}
		\end{enumerate}
Observe that \ref{P:weak} implies that the functions $P,P',P'',\h,\h',\h''$ are \Lip\ continuous on $\erre$.
Let us also note that both $F_1=F_{\rm 1,log}$ and $F_1=\indi$ are admissible for \ref{F:weak}. Moreover, \ref{F:weak} implies that the subdifferential $\partial F_1$ of $F_1$ is a 
maximal monotone graph in $\erre \times\erre$ with effective domain $D(\partial F_1 ) \pier{{}\subseteq{}} D(F_1 ) $;
 since $F_1$ attains its minimum value~$0$ at $0$, it also turns out that 
$0\in D(\partial F_1 )$ and $0\in\partial F_1(0)$.
 
Next, we introduce our notion of solution to the state system \Statesys. 
\begin{definition}
\label{DEF:WEAK}
A \last{quadruplet} $(\m,\ph,\xi,\s)$ {is called} a weak solution to the initial-boundary value problem \Statesys\ if
\begin{align}
	\ph & \in \H1 H \cap \L\infty V \cap \L2 {W_0}, \label{pier2-1}
	\\
	\m,\s & \in \H1 {V^*} \cap \L\infty H \cap \L2 V, \label{pier2-2}
	\\ 
	\xi & \in \L2 H, \label{pier2-3}
\end{align}
and if $(\m,\ph,\xi,\s)$ satisfies the corresponding weak formulation {given by}
\begin{align}
	 \label{var:1}& \<\dt(\alpha\mu + \ph), v > 
	+ \iO \nabla \mu \cdot \nabla v
	= \iO P(\ph)(\sigma+\chi(1-\ph)-\mu)v
	-\iO \h(\ph)u_1 v
	\nonumber\\
	& \qquad \hbox{for every $v \in V $ and \last{almost everywhere} in $(0,T)$,}\\[2mm]
	\label{var:2}&\beta\dt\vp-\Delta\vp+\xi+F_2'(\vp)=\mu+\chi\,\sigma, \quad \hbox{$\xi \in \partial F_1(\ph)$, \, a.e. in $Q$,}
 	\\
 	\label{var:3}
 	& \<	\dt\sigma, v>
 	+ \iO \nabla \sigma\cdot \nabla v
 	=\chi\iO \nabla \vp\cdot \nabla v 
 	- \iO P(\ph)(\sigma+\chi(1-\ph)-\mu)v
 	+ \iO u_2 v	\nonumber\\
	& \qquad \hbox{for every $v \in V $ and \last{almost everywhere} in $(0,T)$,}\\
	\label{var:4}
&	\m(0)=\m_0, \quad
	\ph(0)=\ph_0, \quad
	\s(0)=\s_0 \quad \hbox{a.e. in $\Omega$}.
\end{align}
\end{definition}
Observe that the homogeneous Neumann boundary conditions \eqref{ss4}
are encoded in the condition \eqref{pier2-1} for $\ph$ (by the definition of the space $W_0$) and 
in the variational equalities \eqref{var:1} and \eqref{var:3} for $\mu$ and $\sigma$, 
by the use of the forms $\iO \nabla \mu \cdot \nabla v  $ and $\iO \nabla \sigma \cdot \nabla v $. 
Moreover, let us point out that at this level the control pair $(u_1,u_2)$ just plays the role of two fixed 
forcing terms in \eqref{var:1} and \eqref{var:3}. \last{Let us also mention that the} initial 
conditions~\eqref{var:4} are meaningful 
since  \eqref{pier2-1} and \eqref{pier2-2} ensure that $\ph\in C^0([0,T];V)$ and 
$\mu, \sigma\in C^0([0,T];H)$.

The following result is a special case of \cite[Thm.~2.2]{CSS2}.
\begin{theorem}
\label{THM:WEAK}
Assume that \ref{const:weak}--\ref{P:weak} are fulfilled, let
the initial data satisfy
\begin{align}
	\label{weak:initialdata}
	\m_0, \s_0 \in \Lx2, 
	\quad
	\ph_0 \in \Hx1,
	\quad
	F_1(\ph_0) \in \Lx1,
\end{align}
and suppose that 
\begin{align}
	\label{u:weak}
		(u_1, u_2) \in L^2(Q) \times L^2(Q).
\end{align}
Then there exists at least one solution $(\m,\ph,\xi,\s)$ in the sense of Definition \ref{DEF:WEAK}.
Moreover, if $u_1 \in L^\infty(Q)$ then there is only one such solution.
\end{theorem}

Observe that the above well-posedness result is valid also for the case when $F_1=\indi$. 
This is not the case for the next result  
concerning the existence of strong solutions, which however applies to the 
logarithmic case $F_1=F_{\rm 1,log}$. For this purpose, we consider the following 
additional condition for the nonlinearity $F_1$:

\vspace{2mm}\noindent
\begin{enumerate}[label={\bf (A\arabic{*})}, ref={\bf (A\arabic{*})}]
\setcounter{enumi}{4}
\item \label{F:strong}
\,There exists an interval $\,(r_-,r_+)\,$ with $\,-\infty\le r_-<0<r_+\le +\infty\,$ such that the
\pier{restriction of $F_1$ to $\,(r_-,r_+)\,$ belongs to $\,C^5(r_-,r_+)$ and such that}
\begin{equation}
\label{Fmehr}
\lim_{r\searrow \, \last{ r_-}}F_1'(r)=-\infty, \quad \lim_{r\nearrow \, \last{r_+}} F_1'(r)=+\infty.
\end{equation}
\end{enumerate}

Let us remark that the regularity postulated above for the potential $F_1$ entails that its 
derivative can be defined in the classical manner \pier{in $(r_-,r_+)$}, so  that we will
no longer need to consider a selection $\xi$ in the notion of strong solution below.
Moreover, it will be useful to fix once and for all some $R>0$ such that
\begin{equation}
\label{defUR}
{\cal U}_R:=\left\{\bu =(u_1,u_2)\in L^\infty(Q)^2:\,\|\bu\|_{L^\infty(Q)^2}<R\right\} \supset \uad.
\end{equation}

We then have the following well-posedness result for the state system (where
the equations and conditions are fulfilled almost everywhere in $Q$), which has 
been proved in \cite[Theorem~2.3]{CSS2}:
\Bthm
\label{THM:STRONG}
Suppose that the conditions {\ref{const:weak}}--{\ref{F:strong}} and \eqref{defUR}
are fulfilled, and let the initial data satisfy the conditions
\begin{align}
	\label{strong:initialdata}
&	\m_0,\s_0 \in \Hx1 \cap \Lx\infty, \quad \ph_0 \in { W_0 ,}\\[1mm]
&	r_-<\min_{x\in\overline{\Omega}}\,\vp_0(x) \le
	\max_{x\in\overline{\Omega}}\,\vp_0(x)<r_+.
	\label{strong:sep:initialdata}
\end{align}
Then the state system \eqref{ss1}--\eqref{ss5} has for every $\bu=(u_1,u_2)\in \UR$ a unique strong solution
$(\mu,\vp,\sigma)$ with the regularity
\begin{align}
\label{regmu}
&\mu\in  H^1(0,T;H) \cap C^0([0,T];V) \cap L^2(0,T;W_0)\cap  L^\infty(Q),\\[1mm]
\label{regphi}
&\ph\in W^{1,\infty}(0,T;H)\cap H^1(0,T;V)\cap L^\infty(0,T;W_0) \cap C^0(\overline Q),
\\[1mm]
\label{regsigma}
&\sigma\in H^1(0,T;H)\cap C^0([0,T];V)\cap L^2(0,T;W_0)\cap L^\infty(Q).
\end{align}
Moreover, there is a {constant} $K_1>0$, which depends on $\Omega,T,R,\alpha,\beta$ and the data of the 
system, but not on the choice of $\bu\in \UR $, such that
\begin{align}\label{ssbound1}
&\|\mu\|_{H^1(0,T;H) \cap C^0([0,T];V) \cap L^2(0,T;W_0)\cap L^\infty(Q)}\nonumber\\[1mm]
&+\,\|\ph\|_{W^{1,\infty}(0,T;H)\cap H^1(0,T;V) \cap L^\infty(0,T;W_0)\cap C^0(\overline Q)}
\nonumber\\[1mm]&+\,\|\sigma\|_{H^1(0,T;H) \cap C^0([0,T];V) \cap L^2(0,T;W_0)\cap L^\infty(Q)}\,\le\,K_1\,.
\end{align}
Furthermore, there are constants $r_*,r^*$, which depend on $\Omega,T,R,\alpha,\beta$ and the data of the 
system, but not on the choice of $\bu\in \UR$, such that
\begin{equation}\label{ssbound2}
r_-  <r_*\le\vp(x,t)\le r^*<r_+ \quad\mbox{for all $(x,t)\in \overline Q$}.
\end{equation}
Also, there is some constant $K_2>0$ having the same dependencies as $K_1$ such that
\begin{align}
\label{ssbound3}
	&\max_{i=0,1,2,3}\,\left\|P^{(i)}(\vp)\right\|_{L^\infty(Q)}\,
	+ \max_{i=0,1,2,3}\,\left\|\h^{(i)}(\vp)\right\|_{L^\infty(Q)}\,
	\non\\
	&+\,\max_{i=0,1,2,3,4{,5}}\,\left\|
	F^{(i)}(\vp)\right\|_{L^\infty(Q)} \,\le\,K_2\,.
\end{align}
\Ethm
\Brem
\last{Condition \eqref{ssbound2}, known as the {\em separation property},}
is \an{especially relevant} for the case of singular 
potentials (such as $F_1=F_{\rm 1,log}$). Indeed, it guarantees that the phase variable \last{$\ph$} always stays away
from the critical values $\,r_-,r_+$ that may correspond to the pure phases. \an{Hence}, the singularity \last{of the potential}
is no longer an obstacle for the analysis \an{as} the values of $\ph$ range in some interval in which 
$F_1$ is smooth.

\Erem

Owing to Theorem 2.3, the control-to-state operator $$\, \S:\bu=(u_1,u_2)\mapsto 
(\mu,\ph,\sigma)\,$$ is well defined as a
mapping between ${\cal U}=L^\infty(Q)^2$ and the Banach space specified by the regularity 
results~\eqref{regmu}--\eqref{regsigma}. We now discuss its differentiability properties. The 
results obtained are originally due to \cite{ST} and have been slightly generalized in
\cite{CSS2} to the version reported here. 
For this purpose, some functional analytic preparations are in order. We first define the linear spaces
\begin{align}
\X\,&:=\,X\times \widetilde X\times X, \mbox{\,\,\, where }\nonumber\\
X\,&:=\,H^1(0,T;H)\cap \pier{{}\L\infty V {}}\cap L^2(0,T;W_0)\cap L^\infty(Q), \nonumber\\
\widetilde X\,&:=W^{1,\infty}(0,T;H)\cap H^1(0,T;V)\cap L^\infty(0,T;W_0)\cap C^0(\overline Q),
\label{calX}
\end{align}
which are Banach spaces when endowed with their natural norms. Next, we introduce the linear space
\begin{align}
&\Y\,:=\,\bigl\{(\mu,\ph,\sigma)\in \calX: \,\alpha\dt\mu+\dt\ph-\Delta\mu\in \Liq, 
\,\,\,\beta\dt\ph-\Delta\ph-\mu\in \Liq,\nonumber\\
& \hspace*{14mm}\dt\sigma-\Delta\sigma+\chi\Delta\ph\in\Liq\bigr\},
\label{calY}
\end{align}
which becomes a Banach space when endowed with the norm
\begin{align}
\|(\mu,\ph,\sigma)\|_{\Y}\,:=\,&\|(\mu,\ph,\sigma)\|_\calX
\,+\,\|\alpha\dt\mu+\dt\ph-\Delta\mu\|_{\Liq}
\,+\,\|\beta\dt\ph-\Delta\ph-\mu\|_{\Liq}\nonumber\\
&+\,\|\dt\sigma-\Delta\sigma+\chi\Delta\ph\|_{\Liq}\,.
\label{normY}
\end{align}
Finally, we put 
\begin{align}
&Z:=H^1(0,T;H)\cap L^\infty(0,T;V)\cap L^2(0,T;W_0),
\label{defZ}
\\
&{\cal Z}:= Z\times \widetilde X\times Z.
\label{calZ}
\end{align} 

Now suppose that $\ubar=(\uebar,\uzbar)\in \UR$ is arbitrary and that $(\bmu,\bvp,\bsigma)={\cal S}(\ubar)$.
We then consider the linearization of the state system at $((\uebar,\uzbar),(\bmu,\bvp,\bsigma))$ given by the 
following linear initial-boundary value problem:
\begin{align}\label{aux1}
&\alpha\dt\eta+\dt\rho-\Delta\eta\,=\,P(\bvp)(\zeta-\chi\rho-\eta)+P'(\bvp)(\bsigma+\chi(1-\bvp)-\bmu)\rho
- \h'(\bvp)\,\uebar\,\rho\nonumber\\
&\hspace*{37mm}-\h(\bvp){h_1}\quad\mbox{in }\,Q,\\[1mm]
\label{aux2}
&\beta\dt\rho-\Delta\rho-\eta\,=\,\chi\,\zeta-F''(\bvp)\rho \quad\mbox{ in \,$Q$},
\\[1mm]
\label{aux3}
&\dt\zeta-\Delta\zeta+\chi\Delta\rho\,=\,-P(\bvp)(\zeta-\chi\rho-\eta)-
P'(\bvp)(\bsigma+\chi(1-\bvp)-\bmu)\rho +h_2 \quad\mbox{ in \,$Q$},\\[1mm]
\label{aux4}
&\dn\eta\,=\,\dn\rho\,=\,\dn\zeta\,=\,0 \quad\mbox{ on \,$\Sigma$},\\[1mm]
\label{aux5}
&\eta(0)\,=\,\rho(0)\,=\,\zeta(0)\,=\,0\,\,\,\mbox{ in }\,\Omega\,.
\end{align}
According to \cite[Lem.~4.1]{CSS2} and its proof (see, in particular,
Remark 4.2 and Eqs. (4.37)--(4.39) in \cite{CSS2}), we have the following: 
\begin{align}
&\mbox{The linear system \eqref{aux1}--\eqref{aux5} has for every $\bh=(h_1,h_2)
\in L^2(Q)^2$ a unique solution }
\nonumber\\
&\mbox{$(\eta,\rho,\zeta)\in {\cal Z}$}, \mbox{and the linear mapping 
$\mathbf{h}\mapsto (\eta,\rho,\zeta)$ belongs to ${\cal L}(L^2(Q)^2,{\cal Z})$}.
\label{linsys1}
\\
&\mbox{The linear system \eqref{aux1}--\eqref{aux5} has for every $\mathbf{h}={(h_1,h_2)}\in \Liq^2$ 
a unique solution}\nonumber\\
&\mbox{$(\eta,\rho,\zeta)\in {\mathcal Y}$}, \mbox{and the linear mapping 
$\mathbf{h}\mapsto (\eta,\rho,\zeta)$ belongs to ${\cal L}(\Liq^2,{\cal Y})$}.
\label{linsys2}
\end{align} 
Moreover, we have the following
differentiability result (see \cite[Thm.~4.4]{CSS2}):
\Bthm
\label{THM:FRECHET}
Suppose that the conditions {\ref{const:weak}--\ref{F:strong}} and \eqref{defUR} are fulfilled,
let the initial data $(\m_0,\ph_0,\s_0)$ satisfy \eqref{strong:initialdata} and 
\eqref{strong:sep:initialdata},
and assume that $\overline \bu=(\uebar,\uzbar)\in {\cal U}_R$ is arbitrary and 
$(\bmu,\bvp,\bsigma)={\cal S}(\overline\bu)$. Then the control-to-state 
operator $\S$ is twice continuously Fr\'echet differentiable at $\,\overline\bu\,$ as a mapping from $\,{\cal U}\,$ into 
$\,{\cal Y}$. Moreover, for every ${ \bh=(h_1,h_2)}\in {\cal U}$, the Fr\'echet derivative $\,D \S(\overline\bu)\in 
{\cal L}({\cal U},\Y)\,$ of $\,\S\,$ at
$\,\overline\bu\,$ is given by the identity $\,D\S(\overline\bu){(\bh)}=(\eta,\rho,\zeta)$, where 
$(\eta,\rho,\zeta)$ is the 
unique solution to the linear system \eqref{aux1}--\eqref{aux5}.
\Ethm

\Brem \label{RMK:2.6}
\pier{As} $L^\infty(Q)^2$ is densely embedded in $L^2(Q)^2$,  
the Fr\'echet derivative $\,D \S(\overline \bu)$, which by virtue of the
continuity of the embedding ${\cal Y}\subset {\cal Z}$ also belongs to the space 
${\cal L}(\Liq^2,{\cal Z})$, can be continuously extended to a
linear operator in ${\cal L}(L^2(Q)^2,{\cal Z})$, which we still denote
by $D \S(\overline\bu)$. It then follows from \eqref{linsys1} that also 
for $\bh=(h_1,h_2) \in L^2(Q)^2$ the
identity $\,D\S(\overline\bu){(\bh)}=(\eta,\rho,\zeta)$ is valid. 
\Erem

\Brem
For the explicit form of the second-order Fr\'echet derivative $D^2 \S(\overline\bu)
\in {\cal L}({\cal U},\cal L({\cal U},{\cal Y})),$ we refer the reader
\an{to \cite[Thm.~4.8]{CSS2}}.
\Erem

\section{Deep quench approximation of the state system}
\setcounter{equation}{0}

In this section, we discuss the deep quench approximation
of the state system \Statesys, where we generally assume
that the conditions {\ref{const:weak}}--{\ref{uad}} and \an{\eqref{defUR}--\eqref{strong:initialdata} are fulfilled and that \eqref{strong:sep:initialdata} is satisfied with $(r_-,r_+)=(-1,1)$.} We now consider the state system \Statesys\ for the cases $F_1=\indi$ and  $F_1=\fal$ ($\gamma\in (0,1]$), respectively. 
Since the logarithmic functions $\fal$ satisfy the condition {\ref{F:strong}}, the state 
system \Statesys\ has by Theorem \ref{THM:STRONG} for every $\bu=(u_1,u_2)\in \UR$ and 
$F_1=\fal$, $\gamma\in (0,1]$, a unique solution \last{triplet} $(\mual,\phial,\sial)$ with
the regularity specified by \eqref{regmu}--\eqref{regsigma}. 
By virtue of Theorem \ref{THM:WEAK}, there also exists a unique weak solution 
$(\mu^0,\an{\ph^0}, \xi^0, \sigma^0)$ to the state system \eqref{var:1}--\eqref{var:4}
for $F_1=\indi$ that enjoys the regularity specified by 
\eqref{pier2-1}--\eqref{pier2-3}. 
Clearly, we must have
\begin{equation}
\label{interv}
\pier{-1\le \phial\le 1\,\mbox{ a.e. in \,$Q$,} \, \hbox{ for all $\gamma\in (0,1],$ \, and }
\, -1\le \vp^0\le 1\, \mbox{ a.e. in $\,Q$.}}
\end{equation} 
We introduce the corresponding solution
operators
\begin{align}
&\SAL:\UR\ni \bu\mapsto 
\SAL(\bu)=\big(\SAL^1(\bu),\SAL^2(\bu),\SAL^3(\bu)\big):=(\mual,\phial,\sial) \quad\mbox{for  $\an{\gamma\in (0,1]}$},
\nonumber
\label{defSAL}\\
&\SO:\UR\ni\bu\mapsto \SO(\bu)=\big(\SO^1(\bu),\SO^2(\bu),
\SO^3(\bu),\SO^4(\bu)\big):=(\mu^0,\vp^0,\xi^0,\sigma^0) .
\end{align}

We are now going to investigate the behavior of the family
$\{(\mual,\phial,\sial)\}_{\gamma>0}$ of deep quench approximations for 
$\gamma\searrow0$. \an{We expect that the solution operator $\SAL$ yields an approximation of $\S_0$ as $\gamma \searrow 0.$ This is made rigorous though
the following result.}

\Bthm \label{THM:DQ}
Suppose that the assumptions {\ref{const:weak}}--{\ref{uad}} and \eqref{defUR}--\eqref{strong:sep:initialdata} are fulfilled, and 
let sequences $\{\gamma_n\}\subset (0,1]$ and 
$\{\bu_n\}\subset \uad$ be given such that $\gamma_n\searrow0$ and $\bu_n\to \bu$ weakly-star in ${\cal U}$  as $n\to\infty$ 
for some $\bu\in\uad$. Moreover, let $(\mualn,\phialn,\sialn)=
{\cal S}_{\gamma_n}(\bu_n)$, $n\in\enne$, and $(\mu^0,\vp^0,\xi^0,\sigma^0)=\SO(\bu)$. 
Then, as $n\to\infty$, 
\begin{align}
\label{conmu}
\mualn\to\mu^0&\quad\mbox{weakly-star in }\,X\,\mbox{ and strongly in }\,
C^0([0,T];H),\\
\label{conphi}
\phialn\to \varphi^0&\quad\mbox{weakly-star in \,$\widetilde X \,$ and strongly in }
\,C^0(\overline Q), \\
\label{conxi}
\faln'(\phialn)\to\xi^0&\quad\mbox{weakly-star in }\,L^\infty(0,T;H),\\
\label{consigma}
\sialn \to \sigma^0&\quad\mbox{weakly-star in }\,X\,\mbox{ and strongly in }\,
C^0([0,T];H),
\end{align}
with the denotations introduced in \eqref{calX}.
\Ethm
\Bdim
The sequence $\{\bu_n\}\subset\uad$ forms a bounded subset of ${\cal U}_R$. Now observe that the 
conditions \eqref{strong:initialdata} and \eqref{strong:sep:initialdata} imply
that there is some constant $C_1>0$ such that
\begin{equation}
\label{initialb}
\|\fal(\varphi_0)\|_{C^0(\overline\Omega)}\,+\,\|\fal'(\varphi_0)\|_{
C^0(\overline\Omega)}\,\le\,C_1 \quad\forall\,
\gamma\in (0,1].
\end{equation}
Therefore,  a closer inspection of the a priori estimates carried out  
in the proofs of \an{\cite[Thms.~2.2,~2.3]{CSS2}} reveals that the estimates
\eqref{ssbound1} and \eqref{ssbound3} \pier{(where $F^{(i)} $ are replaced by  $F_2^{(i)} $)} 
hold uniformly for $\gamma\in(0,1]$; in particular, the constant $K_1$ introduced in Theorem~\ref{THM:STRONG} can be chosen in such a 
way that
\begin{align}
\label{albound}
\|\mual\|_X\,+\,\|\phial\|_{\widetilde X} \,+\,\|\sial\|_X\,\le\,K_1
\quad\forall\,\gamma\in (0,1].
\end{align} 
In addition, there is some $C_2>0$ such that
\begin{equation}
\label{Fal}
\|\fal'(\phial)\|_{L^\infty(0,T;H)}\,\le\,C_2\quad\forall \,\gamma\in(0,1].
\end{equation}
Therefore, there are limits $(\mu,\vp,\xi,\sigma)$  and a subsequence of 
$\{(\mualn,\phialn,\sialn)\}$,
which for convenience is again indexed by $n$, such that, as~$n\to\infty$, 
\begin{align}
\label{p1}
\mualn\to\mu&\quad\mbox{weakly-star in }\,X\,\mbox{ and strongly in }\,
C^0([0,T];H),\\
\label{p2}
\phialn\to\varphi&\quad\mbox{weakly-star in }\,\widetilde X\,\mbox{ and strongly in }
\,C^0(\overline Q),\\
\label{p3}
\faln'(\phialn)\to\xi&\quad\mbox{weakly-star in }\,L^\infty(0,T;H),\\
\label{p4}
\sialn \to \sigma&\quad\mbox{weakly-star in }\,X\,\mbox{ and strongly
in }\,\pier{C^0([0,T];H)}.
\end{align}
Here, the strong convergence results follow from well-known compactness results
(see, e.g., \cite[Sect.~8, Cor.~4]{Simon}). 
We \last{then} have to show that $(\mu,\vp,\xi,\sigma)$ is a solution to 
\eqref{var:1}--\eqref{var:4} in the sense of 
Theorem \ref{THM:WEAK} for $F_1=\indi$ and control $\,\bu$. To this end, we pass to the limit as 
$n\to\infty$ in the system \eqref{var:1}--\eqref{var:4},
written for $F_1=\faln$ and $\bu=\bu_n$, for $n\in\enne$. In view of the strong convergence properties stated in \eqref{p1}, 
\eqref{p2}, \last{and} \eqref{p4},
it is easily seen
that $(\mu,\vp,\sigma)$ fulfills the initial \last{conditions in \eqref{var:4}}. Moreover,
owing to the Lipschitz continuity of $P,\h,F_2'$ and the strong convergence 
in \eqref{p2}, we conclude that
\begin{equation}\label{susi}
P(\phialn)\to P(\vp),\quad \h(\phialn)\to \h(\vp),\quad F_2'(\phialn)\to F_2'(\vp),
\quad\mbox{all strongly in }\,C^0(\overline Q).
\end{equation}
Using this and \eqref{p1}--\eqref{p4} \last{once more}, we obtain by passage to the limit as 
$n\to\infty$ that $(\mu,\vp,\xi,\sigma)$ satisfies the time-integrated version
of the variational equalities (with test functions $v\in L^2(0,T;V)$)
stated in \eqref{var:1}--\eqref{var:3}. Notice that this time-integrated version
of the variational equalities is equivalent to them.

It remains to show that $\xi\in\sindi(\vp)$ almost everywhere in $\,Q$. To this end,
we define on $L^2(Q)$ the convex functional
\begin{equation*}
\Phi(v)=\iint_Q \indi(v), \quad\mbox{if \,$\indi(v)\in L^1(Q),\,$ \,\,and 
\,$\Phi(v)=+\infty\,$,  \,otherwise}.
\end{equation*} 
It then suffices to show that $\xi$ belongs to the subdifferential of $\,\Phi\,$
at $\,\vp$, i.e., that
\begin{equation}
\label{reicht}
\Phi(v)-\Phi(\vp)\,\ge\,\iint_Q \xi\,(v-\vp)\quad\forall\,v\in L^2(Q).
\end{equation}
At this point, we recall \eqref{ssbound2} which yields that $\phialn(x,t)\in [-1,1]$
on $\overline Q$.
Hence, by \eqref{p2}, also $\vp(x,t)\in [-1,1]$ on $\overline Q$, and thus 
$\,\Phi(\vp)=0$. Now observe that in case that \,$\Phi(v)\not\in L^1(Q)$ \,the 
inequality \eqref{reicht} holds true since its left-hand side is infinite. If,
however,
$\Phi(v)\in L^1(Q)$, then obviously $v(x,t)\in [-1,1]$
almost everywhere in $\,Q$, and by virtue of \eqref{monoal} and \eqref{limF1al} 
it follows from Lebesgue's dominated convergence theorem that
$$
\lim_{n\to\infty}\iint_Q \faln(v)= \Phi(v)=0.
$$
Now, by the convexity of $\faln$, and since $\faln(\phialn)$ is nonnegative, 
\pier{for all $v\in L^2(Q)$} we have that
$$
\faln'(\phialn)(v-\phialn)\,\le\,\faln(v)-\faln(\phialn)\,\le\,\faln(v)
\quad\mbox{a.e. in \,$Q$.}
$$
Using \eqref{p2} and \eqref{p3}, we thus obtain the following chain of (in)equalities:
\begin{eqnarray*}
\iint_Q \xi(v-\vp) &\!\!=\!\!&\lim_{n\to\infty}\iint_Q \faln'(\phialn)(v-\phialn)
\,\le\,\pier{\limsup_{n\to\infty}}\iint_Q\Bigl(\faln(v)-\faln(\phialn)\Bigr)\\
&\!\!\le\!\!&\lim_{n\to\infty}\iint_Q\faln(v)\,=\,\Phi(v)\,=\,\Phi(v)-\Phi(\vp),
\end{eqnarray*}
which shows the validity of \eqref{reicht}. Hence, the \last{quadruplet}
$(\mu,\vp,\xi,\sigma)$ is a solution to the state system in the sense of
Definition \ref{DEF:WEAK}  for $F_1=\indi$ and the control $\bu$. Since this solution is 
uniquely determined, we must have $(\mu,\vp,\xi,\sigma)=
(\mu^0,\vp^0,\xi^0,\sigma^0)=\SO(\bu)$. Finally, the uniqueness of the limit
also entails that the convergence properties \eqref{p1}--\eqref{p4} are
in fact valid for the entire sequence $(\mualn,\phialn,\sialn)$ and
not only for a subsequence. 
This concludes the proof of the 
assertion.
\Edim
\Brem
Note that the stronger conditions on the data 
\last{required by \eqref{strong:initialdata}--\eqref{strong:sep:initialdata}}
yield more regularity 
for the solution in the case $F_1=\indi$ \pier{with respect to the one} 
obtained from Theorem~\ref{THM:WEAK}. Indeed, we have
\begin{equation}\label{glatter}
\mu\in X, \quad\vp\in \widetilde X,\quad \xi\in L^\infty(0,T;H),
\quad \sigma\in X.
\end{equation}
\Erem

\Brem
\label{Pier}
\pier{The reader may wonder whether a result similar to Theorem~\ref{THM:DQ}
can be proved in the case when the additional assumptions \eqref{strong:initialdata}--\eqref{strong:sep:initialdata}
are not required for the initial data of the state system of \CP0, i.e., of the problem 
\eqref{var:1}--\eqref{var:4} with $F_1 = I_{[-1,1]}$. Indeed, we recall that Theorem~\ref{THM:WEAK} 
states existence and uniqueness of a weak solution to the problem provided the initial data just satisfy \eqref{weak:initialdata}. Note that in this weaker setting the condition  $F_1 (\ph_0) \in L^1(\Omega)$ \last{entails}
\juerg{that} $-1\leq\ph_0\leq 1$ a.e.\ in $\Omega.$ 
The answer to the \juerg{above} question is positive, but in this case the set of initial data $(\m_0,\ph_0,\s_0)$ should be approximated (as $F_1$ is by $F_{1,\gamma}$) by a 
family $\{(\mu_{0,\gamma},\vp_{0,\gamma},\sigma_{0,\gamma})\}$ which does satisfy \eqref{strong:initialdata} and \eqref{strong:sep:initialdata} for every $\gamma \in (0,1]$ and \juerg{converges} to $(\m_0,\ph_0,\s_0)$ in some topology as $\gamma \searrow 0$.
We prove the existence of such a family, with precise statement and all needed conditions, in Lemma~\ref{THM:Pier} in the Appendix.
About the convergence theorem alternative to Theorem~\ref{THM:DQ}, we point out that 
\eqref{conmu}--\eqref{consigma} would hold with the spaces  $X$ and $\widetilde X$ now replaced 
by~(cf.~\eqref{pier2-1}--\eqref{pier2-2})
\begin{align}
&X_w \,:=\,H^1(0,T;V^*)\cap {\L\infty H}\cap L^2(0,T;V), \nonumber\\
&\widetilde X_w \, :=\, \H1H \cap \L\infty V\cap L^2(0,T;W_0)\cap L^\infty (Q),\nonumber
\end{align}
and with $\C0H$ replaced by $\L2H$ in \eqref{conmu} and \eqref{consigma}, $C^0(\overline Q)$ by $\C0H$ in  \eqref{conphi} (and \eqref{susi}), \juerg{and} $\L\infty H $ by $\L2H$ in \eqref{conxi}. Moreover, if one wants to verify the subsequent theory in this weaker setting, it turns out it can be adapted without major modifications 
(see the \last{subsequent} Remark~\ref{Pier2}).}
\Erem 


\section{Existence and approximation of optimal controls}

\setcounter{equation}{0}
Beginning with this section, we investigate the optimal control problem \CP0 of minimizing the cost 
functional \eqref{cost} over the admissible set $\uad$ subject to state system 
\Statesys\ in the form \eqref{var:1}--\eqref{var:4} for $F_1=\indi$. We make the following general assumptions:

\begin{enumerate}[label={\bf (C\arabic{*})}, ref={\bf (C\arabic{*})}]
\item \label{C1} \vspace{2mm}\noindent
  \,\,The constants $b_1,b_2,\kappa$ are nonnegative, and $b_0$ is positive.

\item \label{C2}\vspace{2mm}\noindent
 \,\,It holds $\widehat\vp_\Omega\in\Ldue$ and $\widehat\vp_Q\in L^2(Q)$.

\item \label{C3}\vspace{2mm} \noindent
 \,\,$g:L^2(Q)^2\to\erre$ is nonnegative, continuous and convex on $L^2(Q)^2$.
\end{enumerate}

\vspace{2mm}\noindent
Observe that {\ref{C3}} implies that $\,g\,$ is weakly sequentially lower 
semicontinuous on $L^2(Q)^2$. Moreover, denoting in the following by $\,\partial\,$ 
the subdifferential mapping in $L^2(Q)^2$, it follows from standard convex analysis that 
$\,\partial g\,$ is defined on the entire space $L^2(Q)^2$ and \pier{is} a maximal monotone 
operator. In addition, 
the mapping
$((\mu,\varphi,\sigma),\bu)\mapsto \J((\mu,\varphi,\sigma),\bu)$ 
defined by the cost functional \eqref{cost} is obviously
 continuous and convex 
(and thus weakly sequentially lower semicontinuous) on the
space $\bigl(L^2(Q)\times C^0([0,T];\Ldue)\times L^2(Q)\bigr)
\times L^2(Q)^2$.

In comparison with \CP0, 
we consider for $\gamma>0$ the following control problem:

\vspace{1mm}\noindent
\CPal \,\,\,Minimize $\,{\cal J}((\mu,\vp,\sigma),\bu)\,$
for $\,\bu\in\uad$, subject to  $(\mu,\varphi,\sigma)
=\S_\gamma(\bu)$.

\vspace{1mm}\noindent
We expect that the minimizers  of \CPal\ are for $\gamma\searrow0$ related to minimizers of \CP0.
Prior to giving an affirmative answer to this conjecture, we first show an existence result for \CPal.

\Bprop \label{PROP:EX:OPTCONT:gamma}
Suppose that {\ref{const:weak}}--{\ref{uad}}, {\ref{C1}}--{\ref{C3}}, and \eqref{defUR}--\eqref{strong:sep:initialdata} 
 are satisfied. Then \mbox{{\rm(}${\cal CP}_{\gamma}${\rm)}} has for
every $\gamma\in(0,1]$ a solution.
\Eprop

\Bdim
Let $\gamma\in (0,1]$ be fixed, and assume that a  minimizing sequence 
$\,\{((\mu_n,\vp_n,\sigma_n),\bu_n)\}$ for \CPal\ is given, where $\bu_n\in\uad$ 
and $(\mu_n,\varphi_n,\sigma_n)=\S_\gamma(\bu_n)$ for all $n\in\enne$. 
Since $\{\bu_n\}\subset\uad$, we may without loss of generality assume that
$\,\bu_n\to\bu$ weakly-star in ${\cal U}$ for some $\bu\in\uad$. Moreover, by the
general bound \eqref{ssbound1}, there are a subsequence 
of $\{(\mu_n,\vp_n,\sigma_n)\}\,$ (which is again labeled by
$n\in\enne$) and limit points $\mu,\vp,\sigma$ such that \eqref{p1}, \eqref{p2},
\eqref{p4}, and \eqref{susi} are valid with $(\mualn,\phialn,\sialn)$ replaced
by $(\mu_n,\vp_n,\sigma_n)$. In addition, since $\gamma>0$ is fixed,
we conclude from \eqref{ssbound2} that there are constants 
$r_*(\gamma),r^*(\gamma)$ such that
$$
-1<r_*(\gamma)\le\vp_n\le r^*(\gamma)<1 \quad \mbox{on $\,\overline Q$ \,
for all $n\in\enne$},
$$
from which it also follows that $\,\fal'(\vp_n)\to \fal'(\vp)$ uniformly 
in $\overline Q$ \last{as $n \to \infty$}. We now write the state system \Statesys\ 
for $F_1=\fal$, $(\mu_n,\vp_n,\sigma_n)$, $\bu_n = \juerg{(u_{n,1}, u_{n,2}){}}$, and pass to the limit as
$n\to\infty$, easily arriving at the conclusion that $(\mu,\vp,\sigma)
=\SAL(\bu)$.  Thus, the pair $((\mu,\vp,\sigma),\bu)$ is admissible
for the minimization problem \CPal. The lower semicontinuity properties
of the cost functional then yield that $((\mu,\vp,\sigma),\bu)$ 
is a solution of \CPal.
\Edim

\Bprop \label{PROP:CESARE}
Suppose that {\ref{const:weak}}--{\ref{uad}}, {\ref{C1}}--{\ref{C3}}, and \eqref{defUR}--\eqref{strong:sep:initialdata} 
 are satisfied, and let sequences $\,\{\gamma_n\}\subset (0,1]\,$ and
$\,\{\bu_n\}\subset\uad\,$ be given such \an{that, as $n\to\infty$,} $\,\gamma_n\searrow0\,$ and 
$\,\bu_n\to \bu\,$ weakly-star in ${\cal U}$
for some $\,\bu\in\uad$. Then, 
\begin{align}
\label{cesareuno}
&{\mathcal{J}}(
\SO(\bu),\bu)\,\le\,\liminf_{n\to\infty}\,{\mathcal{J}}
({\cal S}_{\gamma_n}(\bu_n),\bu_n),\\[0.5mm]
\label{cesaredue}
&{\mathcal{J}}(\SO({\mathbf{v}}),\mathbf{v})\,=\,\lim_{n\to\infty}\,
{\mathcal{J}}({\cal S}_{\gamma_n}(\mathbf{v}),\mathbf{v}) \quad\forall\,\mathbf{v}\in\uad.
\end{align}
\Eprop

\Bdim
Theorem \ref{THM:DQ} yields 
that \pier{the component $\phialn$ of 
$\S_{\gamma_n}(\bu_n)=(\mualn,\phialn,\sialn)$ fulfills 
\last{the convergence} \eqref{conphi}.}   
The validity of \eqref{cesareuno} is then a direct consequence of 
the semicontinuity properties of the cost
functional~{\eqref{cost}}.
 
Now suppose that $\mathbf{v}\in\uad$ is arbitrarily chosen, and put 
$(\mualn,\phialn,\sialn):=\S_{\gamma_n}(\mathbf{v})$ for all $n\in\enne$,
 as well as $(\mu^0,\vp^0,\xi^0,\sigma^0):=\SO(\mathbf{v})$.
Applying Theorem \ref{THM:DQ} with the constant sequence $\bu_n=\mathbf{v}$, $n\in\enne$, we see that
\eqref{conmu}--\eqref{consigma} are valid once more. Since the first two summands of the cost functional are continuous with respect to the strong topology of $C^0([0,T];H)$, we
conclude the validity of
\eqref{cesaredue}. 
\Edim

We are now in a position to prove the existence of minimizers for the problem \CP0. We have the following result.

\Bcor \label{COR:EX:OPTCONT:limit}
Suppose that {\ref{const:weak}}--{\ref{P:weak}}, {\ref{C1}}--{\ref{C3}}, and \eqref{defUR}--\eqref{strong:sep:initialdata} 
 are satisfied. Then the optimal control problem {\rm 
($\mathcal{CP}_0$)} has at least
one solution.
\Ecor

\Bdim
Pick an arbitrary sequence $\{\gamma_n\}\subset (0,1]$ such that $\gamma_n\searrow0$ as $n\to\infty$.
By virtue of Proposition \ref{PROP:EX:OPTCONT:gamma}, the optimal control problem (${\mathcal{CP}}_{\gamma_n}$) has for every $n\in\enne$ a solution $((\mualn,\phialn,\sialn),\bu_
{\gamma_n})$ where $(\mualn,\phialn,\sialn)=\S_{\gamma_n}(\bu_{\gamma_n})$ 
for $n\in\enne$.
Since $\uad$ is bounded in ${\cal U}$, we may without loss of generality 
assume that $\bu_{\gamma_n}\to \bu$ weakly-star in ${\cal U}$
for some $\bu\in\uad$.  We then obtain from Theorem \ref{THM:DQ} that \eqref{conmu}--\eqref{consigma}   hold true with $(\mu^0,\vp^0,\xi^0,\sigma^0)=\SO(\bu)$. Invoking the optimality of
$((\mualn,\phialn,\sialn),\bu_{\gamma_n})$ for (${\mathcal{CP}}_\an{\gamma_n}$), we 
then find from Proposition \ref{PROP:CESARE} for every $\,\mathbf{v}\in\uad\,$ the chain of (in)equalities
\begin{align}
&{\cal J}(\SO(\bu),\bu)\,\le\,\liminf_{\an{n}\to\infty}\,
{\cal J}(\S_{\gamma_n}(\bu_{\gamma_n}),\bu_{\gamma_n})\,\le\,
\liminf_{\an{n}\to\infty}\,{\cal J}(\S_{\gamma_n}(\mathbf{v}),\mathbf{v})
\,=\,{\cal J}(\SO(\mathbf{v}),\mathbf{v} ),\non
\end{align}
which yields that $\,(\SO(\bu),\bu)\,$ is an optimal pair for \CP0. The assertion is thus proved.
\Edim

Theorem \ref{THM:DQ} and the proof of Corollary \ref{COR:EX:OPTCONT:limit} indicate that optimal controls of \CPal\ 
are ``close'' to optimal
controls of \CP0 \an{as $\gamma$ approaches zero}. However, they do not yield any information on whether every optimal control
of \CP0\ can be approximated in this way. In fact, such a global result cannot be expected to hold true. 
\an{Nevertheless}, a local answer can be given \an{by employing} a well-known trick. To this end, 
let $\ubar=(\uebar,\uzbar)\in\uad$ be an optimal control for \CP0\ with the associated state 
$\SO(\ubar)$. We associate
with this optimal control the {\em adapted cost functional}
\begin{equation}
\label{adcost}
\widetilde{\cal J}((\mu,\vp,\sigma),\bu):=
{\cal J}((\mu,\vp,\sigma),\bu)\,+\,\frac 12\,\|\bu-\ubar\|^2_{L^2(Q)^2}
\end{equation}
and a corresponding \emph{adapted optimal control problem} for $\gamma>0$, namely:

\vspace{2mm}\noindent
($\widetilde{\mathcal{CP}}_{\gamma}$)\quad Minimize $\,\, 
\widetilde {\cal J}((\mu,\vp,\sigma),\bu)\,\,$
for $\,\bu\in\uad$, subject to the condition that 
$(\mu,\vp,\sigma)=\linebreak \hspace*{14mm}\SAL(\bu)$.

\vspace{2mm}\noindent
With essentially the same proof as that of Proposition \ref{PROP:EX:OPTCONT:gamma} (which needs no repetition here), we can show the 
following result.

\Blem
Suppose that the assumptions of Proposition \ref{PROP:EX:OPTCONT:gamma} are fulfilled. Then the \an{adapted} optimal control problem 
$(\widetilde{\cal CP}_\gamma)$ has for every $\gamma>0$ at least one solution.
\Elem
  
\vspace{1mm}
We are now in the position to give a partial answer to the question raised above \last{through} the following result.

\Bthm \label{THM:APPROX}
Let the assumptions of Proposition \ref{PROP:EX:OPTCONT:gamma} be fulfilled, suppose that 
$\ubar\in \uad$ is an arbitrary optimal control of {\rm $({\mathcal{CP}}_{0})$} with associated state  
$(\bmu,\bvp,\overline\xi,\bsigma)=\S_0(\ubar)$, and let $\,\{\gamma_k\}_{k\in\enne}\subset 
(0,1]\,$ be
any sequence such that $\,\gamma_k\searrow 0\,$ as $\,k\to\infty$. Then there exist a 
subsequence $\{\gamma_{n}\}$ of $\{\gamma_k\}$, and, for every $n\in\enne$, an optimal control
 $\,\bu_{\gamma_{n}}\in \uad\,$ of the adapted problem 
{\rm $(\widetilde{\mathcal{CP}}_{\gamma_{n}})$}
 with associated state $(\mu_{\gamma_{n}},\vp_{\gamma_{n}},\sialn)=\S_{\gamma_{n}}
(\bu_{\gamma_{n}})$,
 such that, as $n\to\infty$,
\begin{align}
\label{tr3.4}
&\bu_{\gamma_{n}}\to \ubar \quad\mbox{strongly in }\,L^2(Q)^2,
\end{align}
and such that \eqref{conmu}--\eqref{consigma} hold true with $(\mu^0,\vp^0,\xi^0,\sigma^0)$ replaced 
by $(\bmu,\bvp,\overline\xi,\bsigma)$. 
 Moreover, we have 
\begin{align}
\label{tr3.5}
&\lim_{n\to\infty}\,\widetilde{{\cal J}}(\S_{\gamma_{n}}(\bu_{\gamma_{n}}),\bu_{\gamma_{n}})
\,=\,{\cal J}(\SO(\ubar),\ubar).
\end{align}
\Ethm

\Bdim
{For any $ k\in\enne$, we pick an optimal control}
$\bu_{\gamma_k} \in \uad\,$ for the adapted problem $(\widetilde{\cal CP}_{\gamma_k})$ and denote by 
$(\mu_{\gamma_k},\vp_{\gamma_k},\sigma_{\gamma_k})=\S_{\gamma_k}(\bu_{\gamma_k})$ the associated strong solution to the 
state system \eqref{ss1}--\eqref{ss5}.
By the boundedness of $\uad$ in $\calU$, there is some subsequence $\{\gamma_{n}\}$ 
of $\{\gamma_k\}$ such that
\begin{equation}
\label{ugam}
\bu_{\gamma_{n}}\to \bu\quad\mbox{weakly-star in }\,{\cal U}
\quad\mbox{as }\,n\to\infty,
\end{equation}
\pier{for} some $\bu\in\uad$. 
Thanks to Theorem \ref{THM:DQ}, the convergence properties \eqref{conmu}--\eqref{consigma} hold true correspondingly
for the \last{quadruple} $(\mu^0,\vp^0,\xi^0,\sigma^0)=\S_0(\bu)$. 
In addition,  the pair $(\S_0(\bu),\bu)$
is admissible for (${\cal CP}_0$). 

We now aim at showing that $ \bu=\ubar$. Once this is shown, it follows from the unique solvability of
the state system \eqref{var:1}--\eqref{var:4} that also $(\mu^0,\vp^0,\xi^0,\sigma^0)=
(\bmu,\bvp, \overline\xi,\bsigma)$. 
Now observe that, owing to the weak sequential lower semicontinuity of 
$\widetilde {\cal J}$, 
and in view of the optimality property of $(\SO(\ubar),\ubar)$ 
for problem $({\cal CP}_0)$,
\begin{align}
\label{tr3.6}
\liminf_{n\to\infty}\, \widetilde{\cal J}(\S_{\gamma_n}(\bu_{\gamma_n}),\bu_{\gamma_n})\,
&\ge \,{\cal J}(\SO(\bu),\bu)\,+\,\frac{1}{2}\,
\|\bu-\ubar\|^2_{L^2(Q)^2}\nonumber\\[1mm]
&\geq \, {\cal J}(\SO(\ubar),\ubar)\,+\,\frac{1}{2}\,\|\bu-\ubar\|^2_{L^2(Q)^2}\,.
\end{align}
On the other hand, the optimality property of  
$\,(\S_{\gamma_{n}}(\bu_{\gamma_n}),\bu_{\gamma_n})
\,$ for problem $(\widetilde {\cal CP}_{\gamma_{n}})$ yields that
for any $n\in\enne$ we have
\begin{equation}
\label{tr3.7}
\widetilde {\cal J}({\cal S}_{\gamma_{n}}(\bu_{\gamma_{n}}),
\bu_{\gamma_{n}})\,\le\,\widetilde {\cal J}({\cal S}_{\gamma_n}
(\ubar),\ubar)\,\pier{{}=\, {\cal J}({\cal S}_{\gamma_n}
(\ubar),\ubar)\,,}
\end{equation}
whence, taking the limit superior as $n\to\infty$ on both sides and invoking (\ref{cesaredue}) in
Proposition \ref{PROP:CESARE},
\begin{align}
\label{tr3.8}
&\limsup_{n\to\infty}\,\widetilde {\cal J}(\S_{\gamma_{n}}(\bu_{\gamma_n}),
\bu_{\gamma_n})\,\le\,\limsup_{n\to\infty}\widetilde {\cal J}(\S_{\gamma_n}(\ubar),\ubar) 
\,=\,\limsup_{n\to\infty} {\cal J}(\S_{\gamma_n}(\ubar),\ubar)
\,=\,{\cal J}(\SO(\ubar),\ubar)\,.
\end{align}
Combining (\ref{tr3.6}) with (\ref{tr3.8}), we have thus shown that 
$\,\frac{1}{2}\,\|\bu-\ubar\|^2_{L^2(Q)^2}=0$\,,
so that $\,\bu=\ubar\,$  and thus also $(\bmu,\bvp,\overline\xi,\bsigma)
=(\mu^0,\vp^0,\xi^0,\sigma^0)$. 
Moreover, (\ref{tr3.6}) and (\ref{tr3.8}) also imply that
\begin{align}
\label{tr3.9}
&{\cal J}(\SO(\ubar),\ubar) \, =\,\widetilde{\cal J}(\SO(\ubar),\ubar)
\,=\,\liminf_{n\to\infty}\, \widetilde{\cal J}(\S_{\gamma_n}(\bu_{\gamma_n}),
 \bu_{\gamma_{n}})\nonumber\\[1mm]
&\,=\,\limsup_{n\to\infty}\,\widetilde{\cal J}(\S_{\gamma_n}(\bu_{\gamma_n}),
\an{\bu_{\gamma_n}}) \,
=\,\lim_{n\to\infty}\, \widetilde{\cal J}(\S_{\gamma_n}(\bu_{\gamma_n}),
\bu_{\gamma_n})\,,
\end{align}                                    
which proves {(\ref{tr3.5})}. 
Moreover, the convergence properties 
\eqref{conmu}--\eqref{consigma} are satisfied. \pier{On the other hand, we have that  
\begin{align}
\label{tr3.9bis}
{\cal J}(\SO(\ubar),\ubar) 
\,&\leq\,\liminf_{n\to\infty}\, {\cal J}(\S_{\gamma_n}(\bu_{\gamma_n}),
 \bu_{\gamma_{n}})
 \,\leq\,\limsup_{n\to\infty}\, {\cal J}(\S_{\gamma_n}(\bu_{\gamma_n}),
 \bu_{\gamma_{n}}) \nonumber\\[1mm]
&\leq \,\limsup_{n\to\infty}\,\widetilde{\cal J}(\S_{\gamma_n}(\bu_{\gamma_n}),
{\bu_{\gamma_n}}) \,
=\,{\cal J}(\SO(\ubar),\ubar) ,
\end{align}        
so that also 
$ {\cal J}(\S_{\gamma_n}(\bu_{\gamma_n}),{\bu_{\gamma_n}})$ converges to $ {\cal J}(\SO(\ubar),\ubar)$ as 
$n\to \infty $, and the relation in \eqref{adcost} enables us to infer \eqref{tr3.4}.}  
\Edim

\section{First-order Necessary Optimality Conditions}
\label{SEC:FOC}
\setcounter{equation}{0}

We now derive first-order necessary optimality conditions for
the control problem \CP0, using the corresponding conditions
for \CPaltil\ as approximations. To this end, we generally assume that the conditions
{\ref{const:weak}}--{\ref{uad}}, {\ref{C1}}--{\ref{C3}}, and \eqref{defUR}--\eqref{strong:sep:initialdata} are fulfilled.
Now let $\ubar\in\uad$ be any fixed optimal control for \CP0\ with
associated state $(\bmu,\bvp,\overline\xi,\bsigma)=\SO(\ubar)$, and assume that $\gamma\in(0,1]$ fixed. 
Moreover, assume that $\ubar_\gamma=\juerg{(\overline u_{\gamma,1},
\overline u_{\gamma,2})}
\in\uad$ is an optimal control for \CPaltil\ with corresponding state
$(\bmu_\gamma,\bvp_\gamma,\bsigma_\gamma)=\SAL(\ubar_\gamma)$.
\pier{Recalling \eqref{cost} \last{and \eqref{adcost}}, we} then
consider the reduced functionals
\begin{align}\label{defG}
&G_1:{\cal U}_R\ni \bu\mapsto \J_1(\SAL(\bu),\bu)+\frac 12 \,\|\bu-\ubar\|^2_{L^2(Q)^2},
\quad G:\UR \ni \bu\mapsto G_1(\bu)+\kappa\,g(\an{\bu})\,.
\end{align}
By Theorem \ref{THM:FRECHET} and the chain rule, $G_1$ is Fr\'echet differentiable at $\ubar_\gamma$, and the Fr\'echet
derivative $DG_1(\ubar_\gamma)\in {\cal L}({\cal U},{\cal Y})$ is given by 
\begin{align}\label{defDG}
DG_1(\ubar_\gamma)(\mathbf{h})=& \,b_1\iint_Q\bigl(\bvp_\gamma-\widehat\varphi_Q
\bigr)\,\rho^{\mathbf{h}}_\gamma
\,+\,b_2\iO\bigl(\bvp_\gamma(T)-\widehat\vp_\Omega\bigr)\,\rho^{\mathbf{h}}_\gamma(T)
\,+\,b_0\iint_Q \ubar_\gamma\cdot\mathbf{h}\nonumber\\
&+\iint_Q(\ubar_\gamma-\ubar)\cdot
\mathbf{h},
\end{align}
for every $\mathbf{h}=(h_1,h_2)$ in ${\cal U}$. Here, the dot
stands for the Euclidean inner
product in $\erre^2$, and 
$(\eta^{\mathbf{h}}_\gamma,\rho^{\mathbf{h}}_\gamma,\zeta^{\mathbf{h}}_\gamma)$
denotes the unique solution  to the linearized system \eqref{aux1}--\eqref{aux5} associated with
$\mathbf{h}=(h_1,h_2)$ and  $(\bmu,\bvp,\bsigma)=(\bmu_\gamma,\bvp_\gamma, \bsigma_\gamma)$.

As in Remark \ref{RMK:2.6}, it follows that $DG_1(\ubar_\gamma)\in{\cal L}({\cal U},{\cal Y})$
can be extended to a linear operator in ${\cal L}(L^2(Q)^2,{\cal Z})$, which is still
denoted by $DG_1(\ubar_\gamma)$ and satisfies \eqref{defDG} for every
$\mathbf{h}=(h_1,h_2)\in L^2(Q)^2$.

Now, by 
arguing along the same lines as in the derivation of \cite[Lem.~3.1]{ST}, we conclude that
there is some $\overline {\boldsymbol \lambda}_\gamma\in \partial g(\ubar_\gamma)
\subset L^2(Q)^2$ such that the following variational inequality is satisfied:
\begin{equation}\label{vug1gamma}
DG_1(\ubar_\gamma)(\mathbf{v}-\ubar_\gamma)\,+\,\kappa\iint_Q \overline{\boldsymbol \lambda}_\gamma\cdot
(\mathbf{v}-\ubar_\gamma)\,\ge\,0\quad\forall\,\mathbf{v}\in\uad.
\end{equation}

As usual, we simplify \eqref{vug1gamma} by means of the adjoint state variables $(p_\gamma,q_\gamma,
r_\gamma)$, which are
defined as the solution \last{triple} $(p,q,r)$ to the adjoint system which is 
formally given by the backward-in-time parabolic system 
\begin{align}
	&\non - \dt p - \b \dt q - \Delta q + \chi \Delta r  + \fal''(\bvp_\gamma)q+F_2''(\bph_\gamma) 
	q + \h'(\bph_\gamma)\,\juerg{\overline u_{\gamma,1}}\, p
	\\ & 
	\label{adj1}
	\quad -P'(\bph_\gamma)(\bs_\gamma+\chi(1-\bph_\gamma)-\bm_\gamma)(p-r) + \chi P(\bph_\gamma)(p-r)
	={{b_1}} (\bph_\gamma - \widehat \ph_Q) &&\mbox{in }\,Q\,,\\[1mm]
	\label{adj2}
	&-\a \dt p-\Delta p -q + P(\bph_\gamma)(p-r)=0 \quad&&\mbox{in }\,Q\,,\\[1mm]
	\label{adj3}
	& - \dt r -\Delta r  - \chi q - P(\bph_\gamma)(p-r) =0 \quad&&\mbox{in }\,Q\,,\\[1mm]
	\label{adj4}
	&\dn p=\dn q=\dn r=0 \quad&&\mbox{on }\,\Sigma\,,\\[1mm] 
	\label{adj5}
	&(p+\b q)(T)= {{b_2}}(\bph_\gamma(T)-\hat \ph_\Omega),
	\quad \a p(T)= 0,\quad r(T)=0\quad &&\mbox{in }\,\Omega\,.
\end{align}
\Accorpa\Adjsys {adj1} {adj5}
Let us point out that the terminal condition for $\,p+\b q\,$ 
prescribe{s} a final datum that only belongs to $\Lx2$.
Therefore, the first equation \eqref{adj1} has to be understood
 in a weak sense.
According \an{to \cite[Thm.~5.2]{CSS2}}, the adjoint system \last{above admits, for every $\gamma$,} a unique 
solution $(\pal,\qal,\ral)$ satisfying
\begin{align}
\label{regsal}
&\pal+\beta\qal\in H^1(0,T;V^*),\\
\label{regpal}
&	\pal  \in \H1 H \cap \L\infty V \cap \L2 {W_0} \cap L^\infty(Q),\\
\label{regqal}
&	\qal \in \L\infty H \cap \L2 V, 	\\
\label{regral}
&	\ral\in \H1 H \cap \L\infty V \cap \L2 {W_0}\cap L^\infty(Q),
\end{align}
such that $(p,q,r)=(\pal,\qal,\ral)$ satisfies
the variational system
\begin{align}
\label{wadj1}
	& -\< \dt (p +\b q),v>
	+ \iO \nabla q \cdot \nabla v
	- \chi \iO \nabla r  \cdot \nabla v
	+ \iO \fal''(\bph_\gamma)\, q\, v +\iO F_2''(\bvp_\gamma)\, q\,v
	\nonumber
	\\
	& \qquad + \iO \h'(\bph_\gamma)\,\ov u_{\gamma,1}\, p\, v
		-\iO P'(\bph_\gamma)(\bs_\gamma+\chi(1-\bph_\gamma)-\bm_\gamma)
		(p-r) v 
	\nonumber\\
	& \qquad + \chi \iO P(\bph_\gamma)(p-r)v = b_1 \iO  (\bph_\gamma - \hat \ph_Q) v,
	\\
	\label{wadj2}
	& -{\a} \iO\dt p\, v
	+ \iO \nabla p \cdot \nabla v 
	-\iO q\,v
	 +\iO  P(\bph_\gamma)\,(p-r)\,v=0,
	\\
	\label{wadj3}
	& 
	- \iO \dt r \,v
	+\iO \nabla r\cdot \nabla v
	- \chi \iO q\, v
	 -\iO  P(\bph_\gamma)\,(p-r)\, v=0 ,
\end{align}
for every $v\in V$ and almost every $t \in (0,T)$, as well as the terminal conditions
\begin{equation}\label{wadj4}
		(p+\b q)(T)=b_2(\bph_\gamma(T)-\hat \ph_\Omega),
		\quad \alpha p(T)= 0,\quad r(T)=0\quad \mbox{a.e. in }\,\Omega\,.
\end{equation}

Now define
\begin{equation}\label{defd}
	{\bf d}_\gamma(x,t) := {\big(\!-\h(\bph_\gamma(x,t)) p_\gamma(x,t) , r_\gamma(x,t)\big)}, \quad 
	\hbox{for a.e. $(x,t) \in Q.$ }
\end{equation}
It is then a standard matter (for the details, see the proof \an{of 
\cite[Thm.~5.4]{CSS2}})
to use the adjoint variables to simplify the variational inequality \eqref{vug1gamma}. It then results
 the following variational inequality:                         
\begin{align}\label{vug2gamma}
	& \iint_Q \bigl({\bf d}_\gamma + b_0\, \ubar_\gamma +\ubar_\gamma -\ubar\bigr)\cdot 
	(\mathbf{v} -\ubar_\gamma)\,+\,\kappa\iint_Q\overline{\boldsymbol \lambda}_\gamma\cdot
	(\mathbf{v}-\ubar_\gamma) 	\geq 	0 \quad \forall\, \mathbf{v} \in \uad,
\end{align}
\an{where $\overline {\boldsymbol \lambda}_\gamma\in \partial g(\ubar_\gamma)\subset L^2(Q)^2$.}
We now pick any sequence $\{\gamma_n\}\subset (0,1]$ such that $\gamma_n\searrow0$. Then, by Theorem \ref{THM:APPROX}, 
we have that (cf. \eqref{conmu}, \eqref{conphi}, \last{and} \eqref{consigma}), \last{as $n\to\infty$,}
\begin{align}\label{ubene}
\ubar_{\gamma_n}\to \ubar &\quad\mbox{strongly in }\,L^2(Q)^2,\\
\label{mubene}
\bmu_{\gamma_n}\to\bmu&\quad\mbox{weakly-star in $\,X$ and strongly in $\,C^0([0,T];L^s(\Omega))\,$
for $\,s\in [1,6)$},\\ 
\label{phibene}
\bvp_{\gamma_n}\to\bvp&\quad\mbox{weakly-star in $\,\widetilde X$ and strongly in }\,C^0(\ov Q),\\
\label{sigmabene}
\bsigma_{\gamma_n}\to\bsigma&\quad\mbox{weakly-star in $\,X$ and strongly in $\,C^0([0,T];L^s(\Omega))\,$
for $\,s\in [1,6)$},
\end{align}
where the strong convergence in \eqref{mubene} and \eqref{sigmabene} follows from 
\cite[Sect.~8,~Cor.~4]{Simon} since $\,V\,$ is compactly embedded in $L^s(\Omega)$ for every  $s\in [1,6)$.
Moreover, we also have, \an{as $n\to\infty,$}
\begin{equation}\label{Tbene}
\bvp_{\gamma_n}(T)\to \bvp(T) \quad\mbox{strongly in }\,C^0(\ov \Omega),
\end{equation}
and, by Lipschitz continuity, \an{that}
\begin{align}
\label{nonlbene}
&F_2''(\bvp_{\gamma_n})\to F_2''(\bvp), \quad P(\bvp_{\gamma_n})\to P(\bvp), \quad P'(\bvp_{\gamma_n})\to P'(\vp),
\quad \h'(\bvp_{\gamma_n})\to \h'(\bvp),\nonumber\\
&\mbox{all strongly in }\,C^0(\ov Q).
\end{align}

We now derive general bounds for the adjoint variables $(\pal,\qal,\ral)$, where we consider the system
\eqref{adj1}--\eqref{adj5} for $(p,q,r)=(\pal,\qal,\ral)$. In this process, we denote by $C_i$, $i\in\enne$,
positive constants that may depend on the data of the system, but not on $\gamma\in(0,1]$. Also, we make 
repeated use of the global \an{(uniform with respect to $\gamma\in(0,1]$)} estimate \eqref{ssbound1} without further reference. We also note that
$-1\le \an{\bvp_\gamma} \le 1$ on $\ov Q$ for all $\pier{\gamma \in (0,1]}$, so that
\begin{align}
	\label{boundal1} 
	& \an{\|F_2''(\bvp_\gamma)\|_{\Liq}\,+\,\max_{i=0,1}\|P^{(i)}(\bvp_\gamma)\|_{\Liq}
	+\|\h'(\bvp_\gamma)\|_{\Liq}
	\le\,C_1\quad\forall\,\gamma\in(0,1].}
\end{align}

\vspace{2mm}\noindent
{\sc First estimate:} 
The following estimate is only formal. For a rigorous proof, it would have to be performed on the level
of a suitable Faedo--Galerkin scheme for the approximate solution of \eqref{adj1}--\eqref{adj5}. For the
sake of brevity, we avoid \an{writing} such a scheme explicitly here and argue formally, knowing that this
estimate can be made rigorous. 

We multiply \eqref{adj1} by $\,\qal$, \eqref{adj2} by $\,-\dt\pal$, \eqref{adj3} by
$\,\chi^2\ral$, add the resulting identities, and integrate over $Q^t:=\Omega\times (t,T)$, where $t\in [0,T)$. 
Then, we add to both sides the same term
$\,\an{\tfrac 12}\|\pal(t)\|^2=\an{\tfrac 12}\|\pal(T)\|^2-
\an{\juerg{\iint}_{Q^t}\pal\,\dt\pal}$. 
Integration by parts then leads to the equality
\begin{align}\label{adj6}
	& \frac {\b}2 \|\qal(t)\|^2 \,+\iint_{Q^t}|\nabla \qal|^2\,+\, \a\iint_{Q^t}|\dt \pal|^2
	+   \frac {1}2 \norma{\pal(t)}_V^2
	+ \frac {\chi^2}2 \norma{\ral(t)}^2\nonumber\\
	&	\quad + \chi^2 \iint_{Q^t}|\nabla \ral|^2 +\iint_{Q^t}\fal''(\bvp_\gamma)\,|\qal|^2\nonumber\\
	& 	=\,
	\frac {\b}2 \norma{\qal(T)}^2 + \pier{{}\frac {1}2 \norma{\pal(T)}_V^2	
	+ \frac {\chi^2}2 \norma{\ral(T)}^2	{}}
	\nonumber\\
	&\quad	
	{}+ b_1\iint_{Q^t} (\bph_\gamma-\hat \ph_Q)\,\qal
	+\chi \iint_{Q^t} \nabla \ral \cdot \nabla \qal
	- \iint_{Q^t} F_2''(\bph_\gamma)\,|\qal|^2 	\nonumber\\ 
	& \quad
	- \iint_{Q^t} \h'(\bph_\gamma)\,\ov u_{\gamma,1} \,\pal\,\qal
	\pier{+ \iint_{Q^t} P'(\bph_\gamma)(\bs_\gamma+\chi(1-\bph_\gamma)-\bm_\gamma)\,(\pal-\ral)\,\qal}\nonumber\\
	&\quad - \chi \iint_{Q^t} P(\bph_\gamma)\,(\pal-\ral)\,\qal
	+\iint_{Q^t} P(\bph_\gamma)\,(\pal-\ral)\,\dt \pal
	- \iint_{Q^t} \pal\, \dt\pal \nonumber\\
	&\quad+ \chi^3 \iint_{Q^t} \qal\,\ral
	+ \chi^2 \iint_{Q^t} P(\bph_\gamma)\,(\pal -\ral)\,\ral
	=: \sum_{i=1}^{\pier{13}}I_i.
\end{align}

\noindent Observe that the last term on the \lhs\ is nonnegative since $\fal''\ge 0$.
We estimate the terms on the \rhs\ individually.
The first \pier{three} terms are bounded by a constant, due to the terminal conditions \eqref{adj5},
the assumption {\ref{C2}}, and the fact that $\,\|\bvp_\gamma\|_{\Liq}\le 1$. 
\last{Likewise, for the fourth term we get
\begin{align*}
	|I_4| \leq C_2 \iint_{Q^t} (|\qal|^2 +1).
\end{align*}
}
Moreover, invoking \eqref{ssbound1},
\eqref{boundal1} and Young's inequality, we easily see that
\begin{align*}
&{|I_5|+|I_6|+|I_7|+|I_8|+|I_9|+|I_{12}|+|I_{13}|}
	\nonumber\\
&\quad \leq \pier{{}\frac {\chi^2}2 \iint_{Q^t} |\nabla\ral|^2 + \frac {1}2 \iint_{Q^t}|\nabla \last{\qal}|^2{}}
+ \last{C_3}\iint_{Q^t}\bigl(|\pal|^2+|\qal|^2+|\ral|^2\bigr) 
\,.
\end{align*}
Finally, owing to \eqref{boundal1} and Young's inequality, 
\begin{align*}
	\pier{|I_{10}|+|I_{11}|} \leq 
		\frac \a 2\iint_{Q^t} |\dt \pal|^2 
	+\last{C_4} \iint_{Q^t}\bigl(|\pal|^2+|\ral|^2\bigr).
\end{align*}
Now, we combine the above estimates and invoke Gronwall's lemma to infer \an{that, for every $\gamma \in (0,1],$}
\begin{align}\label{boundal2}
	\norma{\pal}_{\H1 H \cap \L\infty V}
	+ \norma{\qal}_{\L\infty H \cap \L2 V}
	+ \norma{\ral}_{\L\infty H \cap \L2 V}
	\leq \last{C_5}.
\end{align}

\noindent
{\sc Second estimate:}
We can now rewrite {equation} \eqref{adj3} as a \an{backward-in-time} parabolic equation with \pier{null terminal condition and}  source term 
$f_r{:} = \chi \qal + P(\bph_\gamma)(\pal-\ral)$, which is uniformly bounded in {$\L\infty H$} due to the above estimate. It is then a standard matter to infer that
\begin{align}\label{boundal3}
	\norma{\ral}_{\H1 H \cap \L\infty V \cap \L2 {W_0}}
	\leq \last{C_6} \quad\forall\,\gamma\in(0,1].
\end{align}
In addition, since $\pier{\ral} (T)=0\in L^\infty(\Omega)$, we can apply the regularity result~\cite[Thm.~7.1, p.~181]{LSU} 
to infer that also
\begin{equation}\label{boundal4}
\norma{\ral}_{L^\infty (Q)}	\leq \last{C_7} \quad\forall\,\gamma\in(0,1].
\end{equation}

\noindent
{\sc Third estimate:}
From equation \eqref{adj2} \pier{(see also \eqref{adj5})} and the {parabolic regularity theory}, we similarly recover that
\begin{align}\label{boundal5}
	\norma{\pal}_{\L2 {W_0}{{}\cap L^\infty (Q)}}
	\leq \last{C_8}\quad\forall\,\gamma\in(0,1].
\end{align}

\noindent
{\sc Fourth estimate:}
For the next estimate, we introduce the space 
\begin{equation}
\label{defQ}
{\cal Q}=\{v\in H^1(0,T;H)\cap L^2(0,T;V):v(0)=0\}, 
\end{equation}
which is a closed subspace of $H^1(0,T;H)\cap L^2(0,T;V)$ and thus a Hilbert space. 
Obviously, ${\cal Q}$ is continuously embedded in $C^0([0,T];H)$, and we
have the dense and continuous embeddings ${\cal Q}\subset L^2(0,T;H)\subset {\cal Q}^*$, where it is understood that
\begin{equation}\label{embo}
\langle v,w\rangle_{\cal Q}\,=\,\int_0^T(v(t),w(t))\,dt \quad\mbox{for all \,$w\in{\cal Q}$\, and $\,v\in L^2(0,T;H)$}.
\end{equation}  
Next,
we recall an integration-by-parts formula, which is well known for more regular functions
and was proved in the following form in \cite[Lem.~4.5]{CGSSIAM}:
if $({\cal V},{\cal H},\an{\cal V^*})$ is a Hilbert \last{triple} and 
$$
w\in H^1(0,T;{\cal H})\cap L^2(0,T;{\cal V}) \quad\mbox{and}\quad z\in H^1(0,T;{\cal V}^*)\cap L^2(0,T;{\cal H}),  
$$
then the function $\,t\mapsto (w(t),z(t))_{\cal H}\,$ is absolutely continuous,  and for every $t_1,t_2\in[0,T]$
it holds the formula
\begin{equation}\label{ibparts}
\int_{t_1}^{t_2}\bigl[(\dt w(t),z(t))_{\cal H}+\langle \dt z(t),w(t)\rangle_{\calV}\bigr]\,dt\,=\,
(w(t_2),z(t_2))_{\cal H}-(w(t_1),z(t_1))_{\cal H},
\end{equation}
where $(\,\cdot\,,\,\cdot\,)_{\cal H}$ and $\langle\,\cdot\,,\,\cdot\,\rangle_{\cal V}$ denote the inner product
in ${\cal H}$ and the dual pairing in ${\cal V}$, respectively.

We apply the above result to the special case when ${\cal H}=H$, ${\cal V}=V$, \an{$z=\pal+\beta\qal$, and
$w=v\in {\cal Q}$.}
\pier{Then,} using the terminal condition \eqref{adj5}, \pier{the fact that} \an{$v(0)=0$ by \eqref{defQ}}, \pier{as well as the estimates \eqref{boundal2} and \eqref{ssbound1},
we have that}
\begin{align}
\label{adj8}
&\Big|\int_0^T\langle \dt(\pal+\beta\qal)(t),v(t)\rangle\,dt \, \Big| 
\,\le\,\Big|\iint_Q(\pal+\beta\qal)\,\dt v\Big|\,+\,\bigl|\juerg{(}(\pal+\beta\qal)(T),v(T))\bigr|
\nonumber\\[1mm]
&\quad \le\,\|\pal+\beta\qal\|_{L^2(Q)}\,\|\dt v\|_{L^2(Q)}\,+\,b_2\, \last{\|\bvp_\gamma(T)-\hat\vp_\Omega\|\,\|v(T)\|}
\nonumber\\[1mm]
&\quad \le\,\last{C_9}\,\|v\|_{H^1(0,T;H)}\,+\,\last{C_{10}}\,\|v\|_{C^0([0,T];H)}\,\le\,\last{C_{11}}\,\|v\|_{\cal Q},
\end{align}
which means that
\begin{equation}\label{boundal6}
\|\dt(\pal+\beta\qal)\|_{{\cal Q}^*}\,\le\,\last{C_{12}} \quad\forall\,\gamma\in(0,1].
\end{equation}
At this point, we can conclude from the estimates \eqref{boundal1}, \eqref{boundal2}--\eqref{boundal5}, \eqref{boundal6}, using comparison in \eqref{wadj1}, that the linear mapping 
\begin{equation}\label{defLamal}
\Lambda_\gamma:L^2(Q)\ni v\mapsto \Lambda_\gamma(v)=\an{\iint_Q \fal''(\bvp_\gamma)\,\qal\,v \in \erre}
\end{equation}
satisfies
\begin{equation}\label{boundal7}
\|\Lambda_\gamma\|_{{\cal Q}^*}\,\le\,\last{C_{13}}\quad\forall\,\gamma\in(0,1].
\end{equation}

We now return to the sequence $\gamma_n\searrow0$ introduced above and recall the convergence properties
\eqref{ubene}--\eqref{nonlbene}. Owing to the global estimates \eqref{boundal2}--\eqref{boundal5}, \eqref{boundal7} \pier{and possibly taking another subsequence,}
we may without loss of generality assume that there are limit points $p,q,r,$ \last{and $\Lambda$} such that, as $n\to\infty$,
\begin{align}
\label{conp}
\paln\to p&\quad\mbox{weakly-star in $\,X,$ and strongly in $\,C^0([0,T];L^s(\Omega))\,$ for $\,
s\in [1,6)$},\\
\label{conq}
\qaln\to q&\quad\mbox{weakly-star in }\, L^\infty(0,T;H)\cap L^2(0,T;V),\\
\label{conr}
\raln\to r&\quad\mbox{weakly-star in $\,X,$ and strongly in $\,C^0([0,T];L^s(\Omega))$\, for 
$s\in [1,6)$},\\
\label{conlam}
\Lambda_{\gamma_n}\to \Lambda&\quad\mbox{weakly in }\,{\cal Q}^*,
\end{align}
where \pier{$X$ is defined in \eqref{calX}} and the strong convergence in \eqref{conp} and \eqref{conr} again follows from \cite[Sect.~8,~Cor.~4]{Simon}.

We now perform a passage to the limit as $n\to\infty$ in the adjoint system \eqref{wadj1}--\eqref{wadj4},
written for $\gamma=\gamma_n$ and $(p,q,r)=(\paln,\qaln,\raln)$, for $n\in\enne$. From the convergence results
stated above, it is obvious that, as $n\to\infty$,
\begin{align}
\label{term1}
F_2''(\bvp_{\gamma_n})\,\qaln & \to F_2''(\bvp)\,q &&\mbox{weakly in }\,L^2(Q),\\
\label{term2}
P(\bvp_{\gamma_n})\,(\paln-\raln) &\to P(\bvp)\,(p-r)&&\mbox{weakly in }\,L^2(Q),\\  
\label{term3}
b_1(\bvp_{\gamma_n}-\hat\vp_Q) &\to b_1(\bvp-\hat\vp_Q)&&\mbox{strongly in }\,L^2(Q),\\
\label{term4}
b_2(\bvp_{\gamma_n}(T)-\hat\vp_\Omega) &\to b_2(\bvp(T)-\hat\vp_\Omega)&&\mbox{\pier{strongly} in }\,\Ldue.
\end{align}
A bit less obvious is the fact that also 
\begin{equation}
\label{term5}
\h'(\bvp_{\gamma_n})\,\ov u_{\gamma_n,1}\,\paln\to \h'(\bvp)\,\ov u_1\,p\quad\mbox{weakly in }\,L^2(Q),
\end{equation}
and
\begin{align}
\label{term6}
&P'(\bvp_{\gamma_n})\,(\bsigma_{\gamma_n}+\chi(1-\bvp_{\gamma_n})-\bmu_{\gamma_n})(\paln-\raln)\to
P'(\bvp)\,(\bsigma+\chi(1-\bvp)-\bmu)(p-r)\nonumber\\
&\mbox{weakly in }\,L^2(Q).                 
\end{align}
We only show the validity of \eqref{term5}, since the proof of \eqref{term6} is similar and even simpler. To this end,
recall the strong convergence properties \eqref{ubene} and \eqref{nonlbene}, as well as the fact 
$\paln\to p$ strongly in $C^0([0,T];H)$, in particular. It is then easily verified that for every 
$z\in\Liq$ it holds
$$
\lim_{n\to\infty} \iint_Q \h'(\bvp_{\gamma_n})\,\ov u_{\gamma_n,1}\,\paln\,z\,=\iint_Q \h'(\bvp)\,\ov u_1 p \,z\,,
$$
that is, we have weak convergence in $L^1(Q)$. On the other hand,  
$\,\{\h'(\bvp_{\gamma_n})\,\ov u_{\gamma_n,1}\,\paln\}\,$ is bounded in $L^2(Q)$, whence \eqref{term5} follows.
\pier{%
\Brem
\label{Pier2}
\last{With reference to Remark~\ref{Pier},
let us suggest the reader to read it again; we comment that
\eqref{mubene}--\eqref{sigmabene} should be} replaced by
\begin{align}
\label{mupier}
\bmu_{\gamma_n}\to\bmu&\quad\mbox{weakly-star in $\,X_w$, and strongly in $\,L^2(0,T;L^s(\Omega))\,$
for $\,s\in [1,6)$},\\ 
\label{phipier}
\bvp_{\gamma_n}\to\bvp&\quad\mbox{weakly-star in $\,\widetilde X_w$, and strongly in $\,\C0{L^s(\Omega)}\,$
for $\,s\in [1,6)$},\\
\label{sigmapier}
\bsigma_{\gamma_n}\to\bsigma&\quad\mbox{weakly-star in $\,X_w$, and strongly in $\,L^2(0,T;L^s(\Omega))\,$
for $\,s\in [1,6)$}.
\end{align}
\last{By Lipschitz continuity, it then} turns out that
\begin{align}
\label{nonlpier}
&F_2''(\bvp_{\gamma_n})\to F_2''(\bvp), \quad P(\bvp_{\gamma_n})\to P(\bvp), \quad P'(\bvp_{\gamma_n})\to P'(\vp),
\quad \h'(\bvp_{\gamma_n})\to \h'(\bvp),\nonumber\\
&\mbox{all strongly in $\,\C0{L^s(\Omega)}\,$ for $\,s\in [1,6)$}.
\end{align}
Then, one may directly check that \eqref{term1}--\eqref{term5} still hold, and about \eqref{term6} 
we have that, for instance,  $ P'(\bvp_{\gamma_n})$  converges strongly in $\,\C0{L^4(\Omega)}$, 
$(\bsigma_{\gamma_n}+\chi(1-\bvp_{\gamma_n})\juerg{-\overline\mu_{\gamma_n}})$ 
converges strongly in  $\,L^2(0,T;L^4(\Omega))$, 
$(\paln-\raln)$ converges weakly-star in  $\,L^\infty (0,T;L^4(\Omega))$, and consequently, 
\begin{align}
\label{term6-pier}
&P'(\bvp_{\gamma_n})\,(\bsigma_{\gamma_n}+\chi(1-\bvp_{\gamma_n})-\bmu_{\gamma_n})(\paln-\raln)\to
P'(\bvp)\,(\bsigma+\chi(1-\bvp)-\bmu)(p-r)\nonumber\\
&\mbox{weakly in }\,L^2(0,T;L^{4/3}(\Omega)),                 
\end{align}
with $L^2(0,T;L^{4/3}(\Omega))$ continuously embedded in $\L2 {V^*}$. Then, the limit procedure in the sequel 
can be carried out also in the weaker setting. 
\Erem 
Now,} we apply the integration-by-parts formula \eqref{ibparts} to see that for every $v\in{\calQ}$ it holds that
\begin{align}\label{term7}
&\lim_{n\to\infty}\int_0^T\langle -\dt(\paln+\beta\qaln)(t),v(t)\rangle\,dt\nonumber\\
&\quad =\,\lim_{n\to\infty}\Big(\iint_Q(\paln+\beta\qaln)\,\dt v-
b_2\iO(\bvp_{\gamma_n}(T)-\hat \vp_\Omega)\,v(T)\Big)\nonumber\\
&\quad =\,\iint_Q(p+\beta q)\,\dt v\,-\,b_2\iO(\bvp(T)-\hat\vp_\Omega)\,v(T)\,.
\end{align}   

At this point, we may pass to the limit as $n\to\infty$ to arrive at the following limit system:
\begin{align}
\label{wadj1neu}
&\langle \Lambda,v\rangle_{{\cal Q}}\,=\,-\iint_Q(p+\beta q)\,\dt v\,+\,b_2\iO(\bvp(T)-\hat\vp_\Omega)\,v(T)
	- \iint_Q \nabla q \cdot \nabla v
	\non\\
	&\quad + \chi \iint_Q \nabla r  \cdot \nabla v-\iint_Q F_2''(\bvp)\, q\,	v\,- \iint_Q \h'(\bph)\,\ov u_{1}\, p\, v
		 \nonumber\\
	&\quad +\iint_Q P'(\bph)(\bs+\chi(1-\bph)-\bm)
		(p-r) v
		- \chi \iint_Q P(\bph)(p-r)v 
		\nonumber\\
	&\quad		+ b_1 \iint_Q  (\bph - \hat \ph_Q) v \qquad\mbox{for all \,$v\in{\calQ}$},
	\\[2mm]
	\label{wadj2neu}
	& -{\a} \iO\dt p\,v
	+ \iO \nabla p \cdot \nabla v 
	-\iO q\,v
	 +\iO  P(\bph)\,(p-r)\,v=0 \nonumber\\
	&\quad \qquad\mbox{for all $\,v\in V$, almost everywhere in $(0,T$)},
	\\[2mm]
	\label{wadj3neu}
	& 
	- \iO \dt r\,v
	+\iO \nabla r\cdot \nabla v
	- \chi \iO q\, v
	 -\iO  P(\bph)\,(p-r)\, v=0 \nonumber\\
	&\quad \qquad\mbox{for all $\,v\in V$, almost everywhere in $(0,T)$},\\[2mm]
\label{wadj4neu}
		& \alpha p(T)= 0,\quad r(T)=0\quad \mbox{a.e. in }\,\Omega\,.
\end{align}

Finally, we consider the variational inequality \eqref{vug2gamma} for $\gamma=\gamma_n$, $n\in\enne$.
First observe that the above convergence results certainly imply that 
\begin{align}\label{gauss1}
\mathbf{d}_{\gamma_n}+b_0\ubar_{\gamma_n}+\ubar_{\gamma_n}-\ubar\,\to\,\mathbf{d}+b_0\ubar
\quad\mbox{strongly in }\,L^2(Q)^2,
\end{align}
\an{where $\mathbf{d} \pier{{}: ={}} (-\h(\bvp)p,r)$} \pier{a.e. in $Q$.}

At this point, we recall that the subdifferential $\,\partial g\,$ is defined on the entire space $L^2(Q)^2$
and maximal monotone, and thus a locally bounded operator. Owing to \eqref{ubene}, the sequence
$\{\ov {\boldsymbol \lambda}_{\gamma_n}\}$ introduced \eqref{vug1gamma} is therefore bounded in $L^2(Q)^2$,
and we may without loss of generality assume that there is some $\ov{\boldsymbol \lambda}\in L^2(Q)^2$ such
that $\,\ov {\boldsymbol \lambda}_{\gamma_n}\to \ov{\boldsymbol \lambda}\,$ weakly in $L^2(Q)^2$ \an{as $n\to\infty$}.
A standard argument for maximal monotone operators then yields that $\ov{\boldsymbol \lambda}\in\partial g(\ubar)$,
and passage to the limit as $n\to\infty$ in \eqref{vug2gamma} then shows that the following variational
inequality is satisfied for the limiting variables:
\begin{align}\label{vug3}
& \iint_Q \bigl({\bf d}\,+\, b_0\, \ubar \,+\,\kappa\,\overline{\boldsymbol \lambda})\cdot
	(\mathbf{v}-\ubar) 	\geq 	0 \quad \forall\, \mathbf{v} \in \uad.
\end{align}

Summarizing the above considerations, we have proved the following first-order necessary optimality conditions for \last{the}
optimal control problem \CP0.

\begin{theorem} \label{THM:ADJ:LIMIT}
Suppose that the conditions {\ref{const:weak}}--\pier{{\ref{uad}}}, {\ref{C1}}--{\ref{C3}}, and 
\eqref{defUR}--\eqref{strong:sep:initialdata} are fulfilled, and let $\ubar\in\uad$ be a 
minimizer of the optimal control problem \CP0\
with associate state $(\bmu,\bvp,\ov \xi, \bsigma)=\SO(\ubar)$. Then there exist $p,q,r,\ov {\boldsymbol \lambda},$ \last{and $\Lambda$}
such that the following holds true:\\[1mm]
{\rm (i)} \quad\,\,$p,r\in X$, $q\in L^\infty(0,T;H)\cap L^2(0,T;V)$, $\ov {\boldsymbol \lambda}\in \last{\partial}g(\ubar)$,
$\Lambda\in {\cal Q}^*$.\\[1mm]
{\rm (ii)} \quad The adjoint system \eqref{wadj1neu}--\eqref{wadj4neu} and the variational inequality \eqref{vug3}
are {satisfied}
\linebreak \hspace*{10mm} where $\,{\bf d}=(-\h(\bvp)p,r)$. 
\end{theorem}

\vspace{2mm}
\Brem
\an{(i)} Observe that the adjoint state $(p,q,r)$ and the Lagrange multiplier $\Lambda$ are not unique. However, all possible
choices satisfy \eqref{vug3}.\\[1mm]
\an{(ii)} We have, for every $n\in\enne$, the complementarity slackness condition
\pier{(cf.~\eqref{defLamal})}
$$\pier{\Lambda_{\gamma_n}(\qaln)}=\iint_Q F_{1,\gamma_n}''(\bvp_{\gamma_n})\,|\qaln|^2\,\ge\,0.
$$
Unfortunately, the weak convergence properties of $\{\qaln\}$ do not permit a passage to the limit in this
inequality to
derive a corresponding result for \CP0.
\Erem

\section{Sparsity of optimal controls}
\setcounter{equation}{0}

In this section, we discuss the sparsity of optimal controls, \last{that is}, the possibility that the optimal controls
will vanish in some proper subset of $\,Q$; the form of this subset depends on the particular choice
of the convex function $g$ in the cost functional, while its size depends on the sparsity parameter $\kappa$ \pier{(see~\eqref{cost})}.
We again generally assume that the conditions {\ref{const:weak}}--{\ref{uad}}, {\ref{C1}}--{\ref{C3}}, and
\eqref{defUR}--\eqref{strong:sep:initialdata} are satisfied. Moreover, we assume that
 $\ubar=(\uebar,\uzbar)\in\uad$ is a minimizer of \CP0\ with associated state
$(\bmu,\bvp,\ov\xi,\bsigma)\an{= \S_0(\juerg{\ubar})}$ and adjoint state $(p,q,r)$. Then, the first-order necessary optimality 
condition \eqref{vug3} is satisfied. Since we plan to discuss sparsity properties, we make a further assumption:

\begin{enumerate}[label={\bf (C\arabic{*})}, ref={\bf (C\arabic{*})}]
\setcounter{enumi}{3}
\item \label{C4}\,\,The sparsity parameter $\kappa$ is positive.
\end{enumerate}

\vspace{2mm}\noindent
The sparsity properties will be deduced from the variational inequality \eqref{vug3} and the particular
form of the subdifferential  $\partial g$. In the following argumentation, we closely follow the lines of
\cite[Sect.~4]{ST}; since a detailed discussion is given there, we can afford to be brief here. 

We begin our analysis by introducing  the convex functionals $g$ we are interested in\last{.}

\vspace{2mm}\noindent
{\em Directional sparsity with respect to time:}  Here we use $\,g_T: L^1(0,T;L^2(\Omega)) \to \erre$,
\begin{equation}
\label{gT}
g_T(u) = \|u\|_{L^1(0,T;\Ldue)}= \int_0^T \|u(\cdot,t)\|_{\Ldue}\,dt,
\end{equation}
with the subdifferential \pier{(in $L^2(Q)$, cf.~\ref{C3})}
\begin{equation}\label{dgT}
\partial g_T(u) = \left\{\lambda \in \pier{L^2(Q)} 
: 
\left\{
\begin{array}{ll}
\|\lambda(\cdot,t)\|_{\Ldue} \le 1 &\,\mbox{ if } \, \pier{\|u(\cdot,t)\|_{\Ldue}} = 0\\[1ex]
\lambda(\cdot,t)\,= \displaystyle \frac{u(\cdot,t)}{ \|u(\cdot,t)\|_{\Ldue}}&\,\mbox{ if } \, 
\pier{\|u(\cdot,t)\|_{\Ldue}} \not= 0
\end{array}
\right.
\right\},
\end{equation}
where the properties above are satisfied for almost every $t \in (0,T)$.

\vspace{2mm}\noindent
{\em Directional sparsity with respect to space:}  In this case we use 
$g_\Omega: L^1(\Omega;L^2(0,T)) \to \erre$,
\begin{equation}
\label{gO}
g_\Omega(u) = \|u\|_{\an{L^1(\Omega;L^2(0,T))}}= \iO \|u(x,\cdot)\|_{L^2(0,T)}\,dx,
\end{equation}
with the subdifferential                        
\begin{equation}
\label{dgO}
\partial g_\Omega(u) = \left\{\lambda \in \pier{L^2(Q)}
: 
\left\{
\begin{array}{ll}
\|\lambda(x,\cdot)\|_{L^2(0,T)} \le 1 &\,\mbox{ if } \, \pier{\|u(x,\cdot)\|_{L^2(0,T)}} = 0\\[1ex]
\lambda(x,\cdot)= \displaystyle \frac{u(x,\cdot)}{\|u(x,\cdot)\|_{L^2(0,T)}}&\,\mbox{ if } \,
\pier{\|u(x,\cdot)\|_{L^2(0,T)}} \not= 0
\end{array}
\right.
\right\},
\end{equation}
where the above properties have to be fulfilled for almost every $x \in \Omega$.

\vspace{2mm}\noindent
{\em Spatio-temporal sparsity:} Here we use $g_Q: L^1(Q) \to \erre$, 
\begin{equation}\label{gQ}
g_Q(u) = \|u\|_{L^1(Q)}=\iint_Q |u(x,t)|\, dx\,dt,
\end{equation}
with the subdifferential 
\begin{equation}\label{dgQ}
\partial g_Q(u) = \left\{\lambda \in \pier{L^2(Q)}:\,
\lambda(x,t)  \left\{
\begin{array}{ll}
=1 & \mbox{ if } \, u(x,t) > 0\\
\in [-1,1]& \mbox{ if } \, u(x,t) = 0\\
= -1 & \mbox{ if } \,  u(x,t) < 0\\
\end{array}
\right.
\mbox{, \ a.e. } {(x,t) \in Q}
\right\}.
\end{equation}

\Brem
Observe that in any of the cases $g\in\{g_T,g_\Omega,g_Q\}$ the subdifferential 
\pier{operates} on the entire
space $L^2(Q)$. Moreover, \pier{in the third example it turns out that whenever $\lambda\in\partial g_Q (u)$}, then $|\lambda|\le 1$ almost everywhere in $Q$. 
\Erem

In the following, we concentrate on directional sparsity in time, since this seems to be \an{the} most important
for medical applications; indeed, if an application to medication is considered, directional sparsity in time 
  will allow to 
stop the administration of drugs in certain intervals of time. The subsequent analysis is based on the following
auxiliary sparsity result (see \cite{herzog_stadler_wachsmuth2012,casas_herzog_wachsmuth2017,ST}) for the case
of scalar controls:

\begin{lemma}\label{L4.1} 
Assume that
\begin{equation} \label{C}
C = \{v \in L^\infty(Q): \underline w \le v(x,t) \le {\hat w} \,\,\mbox{ for a.e. $(x,t)$ in } Q\},
\end{equation}
with real numbers $\underline w < 0 < \hat w$, and let a function
$d \in L^2(Q)$ be given. Moreover, assume that $\,u \in C\,$ is a solution to the variational 
inequality
\begin{equation} \label{vugaux}
\iint_Q (d + \kappa \lambda + b_0 u)(v - u) \ge 0 \quad \forall\, v \in C,
\end{equation}
with some $\lambda \in \partial g_T(u)$. 
 Then, for almost every $t\in (0,T)$,
\begin{equation}\label{equivalence}
\pier{\|u(\cdot,t)\|_{L^2(\Omega)}{}}=0 \quad \Longleftrightarrow \quad  \|d(\cdot,t)\|_{L^2(\Omega)} \le \kappa,
\end{equation}
as well as 
\begin{equation} \label{subdiff}
\lambda(\cdot,t) \left\{
\begin{array}{lcl}
\in B(0,1)&\,\mbox{ if } \, \|u(\cdot,t) \|_{L^2(\Omega)} = 0\\[1ex]
= \displaystyle \frac{u(\cdot,t)}{\|u(\cdot,t) \|_{L^2(\Omega)}}&\,\mbox{ if } \, \|u(\cdot,t) \|_{L^2(\Omega)} \not= 0
\end{array}
\right. ,
\end{equation}
where $B(0,1) = \{ v \in L^2(\Omega): \|v\|_{L^2(\Omega)}\le 1\}.$
\end{lemma}

We now apply Lemma \ref{L4.1} to the optimal control problem \CP0\ for which the variational inequality \eqref{vug3}
holds true. To this end, we use the convex continuous functional
\begin{equation} \label{gT2}
g(\bu) =  g(u_1,u_2) : = g_T(u_1) + g_T(u_2) = g_T(I_1\bu) + g_T(I_2\bu),
\end{equation}
where $I_i$ denotes the linear and continuous projection mapping $I_i: \bu=(u_1,u_2)  \mapsto u_i$, $i = 1,2$, 
from $L^1(0,T;L^2(\Omega))^2$ to $L^1(0,T;L^2(\Omega))$.
Since the convex functional $g_T$ is continuous on the whole space $L^1(0,T;L^2(\Omega))$,
we obtain from the well-known sum and chain rules for subdifferentials that 
\[
\partial g(\bu) = I_1^*\, \partial g_T(I_1\bu) + I_2^*\, \partial g_T(I_2\bu) = (I,0)^\top \partial g_T(u_1) + (0,I)^\top \partial g_T(u_2),
\] 
with the identical mapping $I \in \mathcal{L}(L^1(0,T;L^2(\Omega)))$. Therefore, we have
\[
\partial g(\bu) = \{ (\lambda_1,\lambda_2) \in 
\last{L^2(Q)^2}:
\lambda_i \in \partial g_T(u_i), \, i = 1,2\}.
\]

Now observe that the variational inequality \eqref{vug3} is equivalent to two independent variational 
inequalities for $\overline {u}_1$ and $\overline {u}_2$ that have to hold simultaneously, namely,
\begin{eqnarray}
\iint_Q  \left( - \h({\overline \varphi}) \,p+ \kappa \overline \lambda_1 +b_0\, \overline u_1\right)
\left(v - \overline{u}_1\right) &\!\!\ge\!\!& 0 \quad \forall\, v \in C_1, \label{varin1}\\
\iint_Q  \left( r
+ \kappa \overline \lambda_2 +b_0\,\overline u_2\right)\left(v - \overline{u}_2\right) 
&\!\!\ge\!\!& 0 \quad \forall\, v \in C_2, \label{varin2}
\end{eqnarray}
with $\ov {\boldsymbol \lambda}=(\ov\lambda_1,\ov \lambda_2)$,  where the sets $\,C_i$, $i = 1,2$, are defined by
\[
C_i = \{v \in L^\infty(Q): \underline u_i  \le v(x,t) \le \hat u_i \mbox{ for a.a. } (x,t) \in Q\},
\]
and where $\overline \lambda_i $, $i = 1,2$, obey {for almost every $t\in (0,T)$} the conditions
\begin{equation} \label{subdiff2}
\overline \lambda_i(\cdot,t) \left\{
\begin{array}{lcl}
\in B(0,1)&\,\mbox{ if } \, \|\overline u_i(\cdot,t) \|_{L^2(\Omega)} = 0\\[1ex]
= \displaystyle \frac{\overline u_i(\cdot,t)}{\|\overline u_i(\cdot,t) \|_{L^2(\Omega)}}&\,\mbox{ if } \, \|\overline u_i(\cdot,t) \|_{L^2(\Omega)} \not= 0
\end{array}
\right. .
\end{equation}
Applying Lemma \ref{L4.1} to each of the variational inequalities \eqref{varin1} and \eqref{varin2} separately, we arrive at the following result:
\begin{theorem}  {\rm (Directional sparsity in time for \CP0)} \\
\noindent Suppose that the general assumptions 
{\ref{const:weak}}--{\ref{uad}}, {\ref{C1}}--{\ref{C4}} and \eqref{defUR}--\eqref{strong:sep:initialdata} are fulfilled, and assume in addition that the constants $\underline u_i,\widehat u_i $ \an{in \ref{uad}} satisfy
\,$\underline u_i<0<\widehat u_i$, for \,$i=1,2$.
Let $\overline \bu = (\overline u_1,\overline u_2)$ be an optimal control of the problem \CP0 with sparsity functional 
$\,g\,$ defined by \eqref{gT2}, and with associated state $(\overline\mu,\bvp,\overline\xi,\overline\sigma)=\SO(\ubar)$ and adjoint state $(p,q,r)$ having the properties stated
in Theorem \ref{THM:ADJ:LIMIT}. Then there are functions $\overline \lambda_i$, $i=1,2,$ that satisfy \eqref{subdiff2} 
and \eqref{varin1}--\eqref{varin2}.
In addition, for
almost every $t\in (0,T)$, we have that
\begin{eqnarray}
\|\overline u_1(\cdot,t)\|_{L^2(\Omega)} = 0 \quad &\Longleftrightarrow& \quad \|\h(\bvp(\cdot,t))\,
p(\cdot,t)\|_{L^2(\Omega)} \le \kappa, \label{u1sparsity}\\
\|\overline u_2(\cdot,t)\|_{L^2(\Omega)} = 0 \quad &\Longleftrightarrow& \quad \| r(\cdot,t)\|_{L^2(\Omega)} \le \kappa. \label{u2sparsity}
\end{eqnarray}
Moreover, if\, $(p,q,r)$ and $\overline \lambda_1, \overline \lambda_2$ are given, then
the optimal controls
$\overline u_1$,  $\overline u_2$ are for almost every $(x,t)\in Q$ obtained from the \last{pointwise} \an{formulae}
\begin{eqnarray*}
\overline u_1(x,t)& =& \max\left\{\underline u_1, \,\min \left\{\hat u_1, -{b_0}^{-1}\,\left(
-\h(\bvp(x,t))\,p(x,t)+ \kappa \overline\lambda_1(x,t)\right)\right\}\,\right\},\\
\overline u_2(x,t) &=& \max \left\{ \underline u_2 ,\,
\min\left\{ \hat u_2 , -{b_0}^{-1} \left(r(x,t) + \kappa \overline \lambda_2(x,t)\right)\right\}\,\right\}.
\end{eqnarray*}
\end{theorem}
\Brem
By virtue of \eqref{u1sparsity} and \eqref{u2sparsity}, optimal controls may vanish on $\Omega$ for some time intervals. Since the functions 
$t \mapsto \|\h({\overline \varphi}(\cdot,t))
\,p(\cdot,t)\|_{L^2(\Omega)}$
and $t \mapsto \| r(\cdot,t)\|_{L^2(\Omega)}$ are continuous on $[0,T]$, 
it is clear that this is the case at least in all open subintervals where these functions are strictly smaller than $\kappa$. 
We also note that one expects the support of optimal controls to shrink with increasing sparsity parameter $\kappa$, which can hardly be quantified. However, it would be useful to confirm that optimal controls vanish for all sufficiently large values of $\kappa$. Unfortunately, while such a result can be shown for the differentiable approximating 
problems \CPal\ and 
\CPaltil\ (by using an argumentation as in the proof of the corresponding 
\cite[Thm.~4.5]{ST}), we are unfortunately unable to recover the necessary 
uniform bounds for $p$ and $r$ from the adjoint state system 
\eqref{wadj1neu}--\eqref{wadj4neu}.   
\Erem

\vspace{2mm}
We conclude this \juerg{section} by briefly sketching the results 
for the other types of sparsity that are obtained if $g$ is given by  
$g_\Omega\,$ or $g_Q$, respectively. In this respect, we refer to 
\cite[Sect.~4.3]{ST}.

\noindent {\em Spatial sparsity:} With the functional 
\last{$g({\bu}) = g_\Omega(u_1) + g_\Omega(u_2)$,} we may have regions in $\Omega$ where the optimal controls vanish for \last{almost every} $t\in (0,T)$. This is established by 
simply interchanging the roles of $t$ and $x$. For instance, instead of the equivalences \eqref{u1sparsity}, \eqref{u2sparsity},
one obtains for \an{almost every} $x \in \Omega$ that
\begin{eqnarray*}
\|\overline u_1(x,\cdot)\|_{L^2(0,T)} = 0 \quad &\Longleftrightarrow& \quad 
\|\h({\overline \varphi}(x,\cdot))p(x,\cdot)\|_{L^2(0,T)} \le \kappa, \\
\|\overline u_2(x,\cdot)\|_{L^2(0,T)} = 0 \quad &\Longleftrightarrow& \quad 
\| r(x,\cdot)\|_{L^2(0,T)} \le \kappa. 
\end{eqnarray*}

\vspace{2mm}\noindent
{{\em Spatio-temporal sparsity:}  If $g$ is defined 
from $g_Q$ by $g({\mathbf u}) = g_Q(u_1) + g_Q(u_2)$,}
then the equivalence relations 
\begin{eqnarray*}
\overline u_1(x,t) = 0 \quad &\Longleftrightarrow& 
\quad |\h({\overline \varphi}(x,t))\,p(x,t)| 
\le \kappa, \\
\overline u_2(x,t) = 0 \quad &\Longleftrightarrow& \quad |r(x,t)| \le \kappa,
\end{eqnarray*}
can be deduced for almost every $\,(x,t) \in Q$. Therefore, the optimal controls may vanish in certain spatio-temporal subsets of $\,Q$.

\section
{Appendix}

\pier{In this section, we \last{show the approximation result} mentioned in Remark~\ref{Pier}.}

\an{\begin{lemma}
\label{THM:Pier}
Let $(\m_0,\ph_0,\s_0) $ be \pier{a {\last{triple}} of initial data satisfying \eqref{weak:initialdata}
with $F_1 = I_{[-1,1]}$.}
Then, there exists an approximating family $\last{\{(\mu_{0,\gamma},\vp_{0,\gamma},\sigma_{0,\gamma})\}}$, \pier{$\gamma\in (0,1]$,}  which satisfies
\begin{align*}
	\non
	& (\mu_{0,\gamma},\vp_{0,\gamma},\sigma_{0,\gamma}) \in 
		(V \cap L^\infty(\Omega)) \times W_0 \times (V \cap L^\infty(\Omega))
		\quad \forall\, \gamma\in(0,1],
	\\ \non  
	& \hspace{.5cm}|\vp_{0,\gamma}| \leq 1- \tfrac {\gamma}2 \quad \aeO,   \quad \forall\, \gamma\in(0,1],
	\\	 \non
	& (\mu_{0,\gamma},\vp_{0,\gamma},\sigma_{0,\gamma})  \to (\m_0,\ph_0,\s_0) \quad \hbox{strongly in } H \times V \times H
			\quad \hbox{as } \gamma\searrow 0.
\end{align*}
\end{lemma}
\begin{proof}
Let us provide a constructive way to produce a possible family of approximating data.
For the first and third variable $\m_{0,\gamma}$ and $\s_{0,\gamma}$, we proceed in the same fashion
because for both of them we need their boundedness in $V \cap L^\infty(\Omega)$.
For every $\gamma \in (0,1]$, we take as $\m_{0,\gamma}$ and $\s_{0,\gamma}$ the unique solutions to the following elliptic problems:
\begin{align*}
 \begin{cases} 
 	\m_{0,\gamma} - \gamma \Delta \m_{0,\gamma} ={\m}_{0} \,\,  &\hbox{in $\Omega$,}
 	\\
 	\dn \m_{0,\gamma} = 0 \,\,& \hbox{on $\Gamma$},
 \end{cases}
    \quad\quad  
 	\begin{cases} 
 	\s_{0,\gamma} - \gamma \Delta \s_{0,\gamma} ={\s}_{0} \,\,  &\hbox{in $\Omega$,}
 	\\
 	\dn \s_{0,\gamma} = 0 \,\, & \hbox{on $\Gamma$}.
 \end{cases}
\end{align*}
Classical theory entails that $\m_{0,\gamma}$ and $\s_{0,\gamma}$ belong to $W_0 \subset L^\infty(\Omega)$,
 and since \pier{$W_0 \subset V$}, we realize that the approximating data $\m_{0,\gamma}$ and $\s_{0,\gamma}$ do enjoy the required regularity. Moreover, from this construction \pier{it} easily follows that, as $\gamma \searrow 0$,
\begin{align*}
	\m_{0,\gamma} \to \m_0, \quad \s_{0,\gamma} \to \s_0, \quad \hbox{\juerg{both} strongly in $H$.}
\end{align*}
We are then reduced to suggest how to construct $\ph_{0,\gamma}$. It is worth recalling that for this \juerg{matter}
 we also 
need to fulfill \eqref{strong:sep:initialdata} with $r_{\pm}= \pm 1$ so that we require the absolute value of $\ph_{0,\gamma}$ to be bounded by $1-\frac {\gamma}2$.
In this direction, we \juerg{introduce} $\hat{\ph}_{0,\gamma}$, the $(1-\frac {\gamma}2)$-truncation of $\ph_0$, that is,
\begin{align*}
	\hat \ph_{0,\gamma} := 
	\begin{cases} 
				 1-\frac {\gamma}2 &\hbox{ if } \ph_0  > {1}-\frac {\gamma}2
				\\ \ph_{0,\gamma}  &\hbox{ if } |\ph_{0,\gamma}| \leq {1}-\frac {\gamma}2
				\\ -{1} +\frac {\gamma}2  &\hbox{ if } \ph_0  < -{1}+\frac {\gamma}2
    \end{cases}.
\end{align*}  
\juerg{Observe that $\,\hat \ph_{0,\gamma}\in V$.}
Then, for every $\gamma \in (0,1]$, we define $\ph_{0,\gamma}$ as the unique solution to the following elliptic problem:
\begin{align}
 \begin{cases} 
 	\ph_{0,\gamma} - \gamma \Delta \ph_{0,\gamma} =\hat{\ph}_{0,\gamma} \quad  &\hbox{in $\Omega$,}
 	\\
 	\dn \ph_{0,\gamma} = 0\quad & \hbox{on $\Gamma$},
 	\label{singular:pert}
 \end{cases}
\end{align}
whose weak formulation is given by
\begin{align}\label{wk:pert}
	\iO \ph_{0,\gamma} \,v +  \gamma \iO \nabla  \ph_{0,\gamma} \cdot \nabla v = \iO \hat{\ph}_{0,\gamma}  \,v \quad \hbox{for every $v \in V $ and a.e. in $(0,T)$}.
\end{align}
Classical theory ensures that, for every $\gamma\in (0,1]$, there exists a unique solution $\ph_{0,\gamma} \in W_0
\juerg{{}\cap H^3(\Omega){}}$.
Moreover, by adding $1-\frac {\gamma}2$ to both sides of the first equation of \eqref{singular:pert}, we 
arrive at the identity
\begin{align*}
	\big(\ph_{0,\gamma} +1-\tfrac {\gamma}2\big) - \gamma \Delta \big(\ph_{0,\gamma} +1-\tfrac {\gamma}2\big)=\hat{\ph}_{0,\gamma} +1-\tfrac {\gamma}2.
\end{align*}
Besides, by construction we have that the \rhs\ is nonnegative so that the maximum principle yields that 
\begin{align*}
	\ph_{0,\gamma} +1-\tfrac {\gamma}2 \geq 0,
	\quad \hbox{ which entails }\quad 
	\ph_{0,\gamma} \geq -1+ \tfrac {\gamma}2. 
\end{align*}
The same strategy can be applied to show also that $\ph_{0,\gamma} \leq 1 - \tfrac {\gamma}2$.
Lastly, it \pier{is} enough to check that $\ph_{0,\gamma}  \to \ph_0$ strongly in $V$ as $\gamma \searrow 0$.
We then test \eqref{wk:pert} by $\ph_{0,\gamma}- \Delta \ph_{0,\gamma}\juerg{{}\in V}$, \juerg{integrate by parts, and} use Young's inequality to find that
\begin{align*}
	& \iO |\ph_{0,\gamma}|^2
	+\iO |\nabla \ph_{0,\gamma}|^2
	+ \gamma \iO |\nabla \ph_{0,\gamma}|^2
	+ \gamma \iO |\Delta \ph_{0,\gamma}|^2
	\\ & \quad \leq 
	\frac 12 \iO |\ph_{0,\gamma}|^2
	+ \frac 12 \iO |\hat{\ph}_{0,\gamma}|^2
	+\frac 12 \iO |\nabla \ph_{0,\gamma}|^2
	+ \frac 12 \iO |\nabla\hat{\ph}_{0,\gamma}|^2,
\end{align*}
whence, \pier{rearranging and multiplying by $2$, we have} that
\begin{align*}
	& \norma{\ph_{0,\gamma}}_V^2
	+ \gamma \iO |\nabla \ph_{0,\gamma}|^2
	+ \gamma \iO |\Delta \ph_{0,\gamma}|^2
	 \leq \norma{\hat{\ph}_{0,\gamma}}_V^2 
	 \leq \norma{\ph_0}_V^2.
\end{align*}
Thus, up to a possible subsequence, we infer that $\ph_{0,\gamma}  \to \ph_0$ weakly in $V$ as $\gamma \searrow 0$, where the identification of the limit easily follows from passing to the limit \pier{as} $\gamma\searrow 0$ in \eqref{wk:pert}.
Furthermore, we infer that $\norma{\ph_{0,\gamma}}^2_V \leq \norma{{\ph}_{0}}^2_V$,
and since this inequality is preserved when passing to the superior limit, we conclude that $\ph_{0,\gamma}  \to \ph_0$ strongly in $V$ as $\gamma\searrow 0$,
as claimed.
\end{proof}}


\section*{Acknowledgments}
\pier{This research was supported by the Italian Ministry of Education, 
University and Research~(MIUR): Dipartimenti di Eccellenza Program (2018--2022) 
-- Dept.~of Mathematics ``F.~Casorati'', University of Pavia. 
In addition, {PC and AS gratefully mention} their affiliation
to the GNAMPA (Gruppo Nazionale per l'Analisi Matematica, 
la Probabilit\`a e le loro Applicazioni) of INdAM (Isti\-tuto 
Nazionale di Alta Matematica). Moreover, PC aims to point out his collaboration,
as Research Associate, to the IMATI -- C.N.R. Pavia, Italy.}

\End{document}